\documentclass{article}

\usepackage{amsfonts,latexsym,amstext}
\usepackage{amsmath,amscd,amssymb}

\usepackage[usenames]{color}

\usepackage[a4paper]{geometry}

\usepackage{url}

\usepackage{psfrag}

\usepackage{placeins}

\usepackage[dvips]{graphicx}

\usepackage{latexsym,amssymb,amsmath,color,times,bbm}

\newtheorem{thm}{Theorem}[section]
\newtheorem{cor}[thm]{Corollary}
\newtheorem{lem}[thm]{Lemma}
\newtheorem{prop}[thm]{Proposition}
\newtheorem{rem}[thm]{Remark}

\numberwithin{equation}{section}

\newcommand{\abs}[1]{\left\vert#1\right\vert}
\newcommand{\set}[1]{\left\{#1\right\}}

\newcommand{\eps}{\varepsilon}

\DeclareMathOperator{\trace}{Tr}

\def\R{\mathbb R}

\def\cprime{$'$}

\newcommand{\tr}[1]{{\vphantom{#1}}^{\mathit t}{#1}}

\def\sqw{\hbox{\rlap{\leavevmode\raise.3ex\hbox{$\sqcap$}}$%
\sqcup$}}
\def\sqb{\hbox{\hskip5pt\vrule width4pt height6pt depth1.5pt%
\hskip1pt}}

\def\qed{\ifmmode\hbox{\hfill\sqb}\else{\ifhmode\unskip\fi%
\nobreak\hfil
\penalty50\hskip1em\null\nobreak\hfil\sqb
\parfillskip=0pt\finalhyphendemerits=0\endgraf}\fi}
\def\cqfd{\ifmmode\sqw\else{\ifhmode\unskip\fi\nobreak\hfil
\penalty50\hskip1em\null\nobreak\hfil\sqw
\parfillskip=0pt\finalhyphendemerits=0\endgraf}\fi}

\begin{document}

\renewcommand{\labelitemi}{$\bullet$}
\bibliographystyle{plain}
\pagestyle{headings}
\title{Ergodic BSDEs and related PDEs with Neumann boundary conditions}
\author{
Adrien Richou\\
IRMAR, Universit\'e Rennes 1\\ Campus de Beaulieu, 35042 RENNES
Cedex, France\\ e-mail: adrien.richou@univ-rennes1.fr}

\maketitle
\abstract{We study a new class of ergodic backward stochastic differential equations (EBSDEs for short) which is linked with semi-linear Neumann type boundary value problems related to ergodic phenomenas. The particularity of these problems is that the ergodic constant appears in Neumann boundary conditions. We study the existence and uniqueness of solutions to EBSDEs and the link with partial differential equations. Then we apply these results to optimal ergodic control problems.}

\section{Introduction}
In this paper we study the following type of (Markovian) backward stochastic differential equations with infinite horizon that we shall call ergodic BSDEs or EBSDEs for short: for all $0 \leqslant t \leqslant T < +\infty,$
\begin{equation}
\label{intro EDSRE gene}
Y^x_t=Y^x_T+ \int_t^T[\psi(X^x_s,Z^x_s)-\lambda]ds+\int_t^T[g(X^x_s)-\mu]dK^x_s-\int_t^T Z^x_s dW_s.
\end{equation}
In this equation $(W_t)_{t\geqslant 0}$ is a $d$-dimensional Brownian motion and $(X^x,K^x)$ is the solution to the following forward stochastic differential equation reflected in a smooth bounded domain $G=\set{\phi>0}$, starting at $x$ and with values in $\mathbb{R}^d$:
\begin{equation}
\label{intro EDS reflechie}
\begin{array}{l}
X^x_t=x+\int_0^t b(X^x_s)ds+\int_0^t \sigma(X^x_s)dW_s+\int_0^t \nabla \phi(X^x_s) dK_s^x, \quad t \geqslant 0;\\
K_t^x = \int_0^t \mathbbm{1}_{X^x_s \in \partial G} dK_s^x, \quad K^x \textrm{ is non decreasing.}
\end{array}
\end{equation}
Our aim is to find a triple $(Y,Z,\mu)$, where $Y,Z$ are adapted processes taking values in $\mathbb{R}$ and $\mathbb{R}^{1 \times d}$ respectively. $\psi: \mathbb{R}^d \times \mathbb{R}^{1 \times d} \rightarrow \mathbb{R}$ is a given function. Finally, $\lambda$ and $\mu$ are constants: $\mu$, which is called the ``boundary ergodic cost'', is part of the unknowns while $\lambda$ is a given constant.
\paragraph{}
It is now well known that BSDEs provide an efficient alternative tool to study optimal control problems, see, e.g. \cite{Peng-93} or \cite{ElKaroui-Mazliak-97}. But up to our best knowledge, the paper of Fuhrman, Hu and Tessitore \cite{Fuhrman-Hu-Tessitore-07} is the only one in which BSDE techniques are applied to optimal control problems with ergodic cost functionals that are functionals depending only on the asymptotic behavior of the state (see e.g. costs defined in formulas (\ref{intro cout I}) and (\ref{intro cout J}) below). That paper deals with the same type of EBSDE as equation (\ref{intro EDSRE gene}) but without boundary condition (and in infinite dimension): their aim is to find a triple $(Y,Z,\lambda)$ such that for all $0 \leqslant t \leqslant T < +\infty,$
\begin{equation}
\label{intro condition de neumann nulle}
Y^x_t=Y^x_T+ \int_t^T[\psi(X^x_s,Z^x_s)-\lambda]ds-\int_t^T Z^x_s dW_s,
\end{equation}
where $(W_t)_{t\geqslant 0}$ is a cylindrical Wiener process in a Hilbert space and $X^x$ is the solution to a forward stochastic differential equation starting at $x$ and with values in a Banach space. In this case, $\lambda$ is the ``ergodic cost''.

There is a fairly large amount of literature dealing by analytic techniques with optimal ergodic control problems without boundary conditions for finite dimensional stochastic state equations. We just mention papers of Arisawa and Lions \cite{Arisawa-Lions-98} and Arisawa \cite{Arisawa-97}. In this framework, the problem is treated through the study of the corresponding Hamilton-Jacobi-Bellman equation. Of course, same questions have been studied in bounded (or unbounded) domains with suitable boundary conditions. For example we refer the reader to Bensoussan and Frehse \cite{Bensoussan-Frehse-02} in the case of homogeneous Neumann boundary conditions and to Lasry and Lions \cite{Lasry-Lions-89} for state-constraint boundary conditions. But in all these works, the constant $\mu$ does not appear and the authors are interested in the constant $\lambda$ instead.

To the best of our knowledge, only works where the problem of the constant $\mu$ appears in the boundary condition of a bounded domain are those of Arisawa \cite{Arisawa-03} and Barles and Da Lio \cite{Barles-DaLio-05}. The purpose of the present paper is to show that backward stochastic differential equations are an alternative tool to treat such ``boundary ergodic control problems''. It is worth pointing out that the role of the two constants are different: our main results say that, for any $\lambda$ and under appropriate hypothesis, there exists a constant $\mu$ for which (\ref{intro EDSRE gene}) has a solution. At first sight $\lambda$ doesn't seem to be important and could be incorporated to $\psi$, but our proof strategy needs it: we first show that, for any $\mu$, there exists a unique constant $\lambda:=\lambda(\mu)$ for which (\ref{intro EDSRE gene}) has a solution and then we prove that $\lambda(\mathbb{R})=\mathbb{R}$.

\paragraph{}
To be more precise, we begin to deal with EBSDEs with zero Neumann boundary condition in a bounded convex smooth domain. As in \cite{Fuhrman-Hu-Tessitore-07}, we introduce the class of strictly monotonic backward stochastic differential equations
\begin{equation}
Y_{t}^{x,\alpha}=Y_{T}^{x,\alpha}+ \int_t^T[\psi(X^x_s,Z_{s}^{x,\alpha})-\alpha Y_{s}^{x,\alpha}]ds-\int_t^T Z_{s}^{x,\alpha} dW_s, \quad 0 \leqslant t \leqslant T < +\infty,
\end{equation}
with $\alpha>0$ (see \cite{Briand-Hu-98} or \cite{Royer-04}). We then prove that, roughly speaking, $(Y^{x,\alpha}-Y^{0,\alpha}_0,Z^{x,\alpha},\alpha Y^{0,\alpha}_0)$ converge, as $\alpha \rightarrow 0$, to a solution $(Y^x,Z^x,\lambda)$ of EBSDE (\ref{intro condition de neumann nulle}) for all $x \in G$ when $(X^x,K^x)$ is the solution of (\ref{intro EDS reflechie}) (see Theorem~\ref{existence convexe}). When there is non zero Neumann boundary condition, we consider a function $\tilde{v}$ such that $\frac{\partial \tilde{v}}{\partial n}(x)+g(x)=\mu,  \forall x \in \partial G$ and thanks to the process $\tilde{v}(X^x)$ we modify EBSDE (\ref{intro EDSRE gene}) in order to apply previous results relating to zero Neumann boundary condition. In Theorems~\ref{existence solution gdiff0} and \ref{unicite lambda gdiff0} we obtain that for any $\mu$, there exists a unique constant $\lambda:=\lambda(\mu)$ for which (\ref{intro EDSRE gene}) has a solution. $\mu \mapsto \lambda(\mu)$ is a continuous decreasing function and, under appropriate hypothesis, we can show that $\lambda(\mu) \stackrel{\mu \rightarrow +\infty}{\longrightarrow} - \infty $ and $\lambda(\mu) \stackrel{\mu \rightarrow -\infty}{\longrightarrow} + \infty$ which allow us to conclude: see Theorem~\ref{existence solution gdiff0 mufixe} when $\psi$ is bounded and Theorems~\ref{existence solution gdiff0 mufixe 2} and \ref{existence solution gdiff0 mufixe 3} when $\psi$ is bounded in $x$ and Lipschitz in $z$. All these results are obtained for a bounded convex domain but it is possible to prove some additional results when the domain is not convex.

Moreover we show that we can find a solution of (\ref{intro EDSRE gene}) such that $Y^x=v(X^x)$ where $v$ is Lipschitz and is a viscosity solution of the elliptic partial differential equation (PDE for short)
\begin{equation}
\label{intro EDP}
\left\{
\begin{array}{l}
\mathcal{L}v(x)+\psi(x,^t\nabla v(x)\sigma(x))=\lambda, \quad x \in G\\
\frac{\partial v}{\partial n}(x)+g(x)=\mu, \quad x \in \partial G,
\end{array}
\right.
\end{equation}
with
\begin{equation*}
\mathcal{L}f(x)=\dfrac{1}{2} \textrm{Tr}(\sigma(x)^t\sigma(x)\nabla^2f(x))+^tb(x)\nabla f(x).
\end{equation*}

The above results are then applied to control problems with costs
\begin{equation}
 \label{intro cout I}
I(x,\rho)=\limsup_{T \rightarrow +\infty} \frac{1}{T} \mathbb{E}^{\rho,T} \left[ \int_0^T L(X_s^x,\rho_s)ds+\int_0^T[g(X_s^x)-\mu]dK_s^x \right],
\end{equation}
\begin{equation}
 \label{intro cout J}
J(x,\rho)=\limsup_{T \rightarrow +\infty} \frac{1}{\mathbb{E}^{\rho,T}[K_T^x]} \mathbb{E}^{\rho,T} \left[ \int_0^T [L(X_s^x,\rho_s)-\lambda]ds+\int_0^Tg(X_s^x)dK_s^x \right] \mathbbm{1}_{\mathbb{E}^{\rho,T}[K_T^x]>0},
\end{equation}
where $\rho$ is an adapted process with values in a separable metric space $U$ and $\mathbb{E}^{\rho,T}$ denotes expectation with respect to $\mathbb{P}_T^{\rho}$ the probability under which $W_t^{\rho}=W_t+\int_0^tR(\rho_s)ds$ is a Wiener process on  $[0,T]$. $R: U \rightarrow \mathbb{R}^d$ is a bounded function. With appropriate hypothesis and by setting  $\psi(x,z)=\inf_{u \in U}\set{L(x,u)+zR(u)}$ in (\ref{intro EDSRE gene}) we prove that $\lambda=\inf_{\rho} I(x,\rho)$ and $\mu=\inf_{\rho} J(x,\rho)$ where the infimum is over all admissible controls.

\paragraph{}
The paper is organized as follows. In the following section we study EBSDEs with zero Neumann boundary condition. In section 3 we treat the general case of EBSDEs with Neumann boundary condition. In section 4 we study the example of reflected Kolmogorov processes for the forward equation. In section 5 we examine the link between our results on EBSDEs and solutions of elliptic semi-linear PDEs with linear Neumann boundary condition. Section 6 is devoted to optimal ergodic control problems and the last section contains some additional results about EBSDEs on a non-convex bounded set.

\section[EBSDEs with zero Neumann boundary condition]{Ergodic BSDEs (EBSDEs) with zero Neumann boundary conditions}
Let us first introduce some notations. Throughout this paper, $(W_t)_{t\geqslant 0}$ will denote a $d$-dimensional Brownian motion, defined on a probability space $(\Omega,\mathcal{F},\mathbb{P})$. For $t\geqslant 0$, let $\mathcal{F}_t$ denote the $\sigma$-algebra $\sigma(W_s;0\leqslant s\leqslant t)$, augmented with the $\mathbb{P}$-null sets of $\mathcal{F}$. The Euclidean norm on $\mathbb{R}^d$ will be denoted by $|.|$. The operator norm induced by $|.|$ on the space of linear operator is also denoted $|.|$. Given a function $f: \mathbb{R}^d \rightarrow \mathbb{R}^k$ we denote $|f|_{\infty}=\sup_{x \in \mathbb{R}^d} |f(x)|$ and $|f|_{\infty,\mathcal{O}}=\sup_{x \in \mathcal{O}} |f(x)|$ with $\mathcal{O}$ a subset of $\mathbb{R}^d$.\\
Let $\mathcal{O}$ be an open connected subset of $\mathbb{R}^d$. $\mathcal{C}^k(\overline{\mathcal{O}})$, $\mathcal{C}^k_b(\overline{\mathcal{O}})$ and $\mathcal{C}^k_{lip}(\overline{\mathcal{O}})$ will denote respectively the set of real functions of class $\mathcal{C}^k$ on $\overline{\mathcal{O}}$, the set of the functions of class $\mathcal{C}^k$ which are bounded and whose partial derivatives of order less than or equal to $k$ are bounded, and the set of the functions of class $\mathcal{C}^k$ whose partial derivatives of order $k$ are Lipschitz functions.\\
$\mathcal{M}^2(\mathbb{R}^+,\mathbb{R}^k)$ denotes the space consisting of all progressively measurable processes $X$, with values in $\mathbb{R}^k$ such that, for all $T>0$,
$$\mathbb{E}\left[ \int_0^T |X_s|^2 ds\right]<+\infty.$$

Throughout this paper we consider EBSDEs where forward equations are stochastic differential equations (SDEs for short) reflected in a bounded subset $G$ of $\R^d$. To state our results, we use the following assumptions on $G$:
\paragraph{(G1).} There exists a function $\phi \in \mathcal{C}_b^2(\mathbb{R}^d)$ such that $G=\set{\phi >0}$, $\partial G=\set{\phi=0}$ and $|\nabla\phi (x)|=1$, $\forall x \in \partial G$.
\paragraph{(G2).} $G$ is a bounded convex set.
\paragraph{} If $x \in \partial G$, we recall that $-\nabla \phi (x)$ is the outward unit vector to $\partial G$ in $x$. We also consider $b: \mathbb{R}^d \mapsto \mathbb{R}^d$ and $\sigma: \mathbb{R}^d \mapsto \mathbb{R}^{d \times d}$, two functions verifying classical assumptions:
\paragraph{(H1).}
there exist two constants $K_b>0$ and $K_{\sigma}>0$ such that $\forall x,y \in \mathbb{R}^d$,
\begin{displaymath}
\begin{array}{crcl}
&|b(x)-b(y)|&\leqslant & K_b |x-y|,\\
\textrm{and} & & &\\
&|\sigma(x)-\sigma(y)|&\leqslant & K_{\sigma} |x-y|.\\
\end{array}
\end{displaymath}
We can state the following result, see e.g. \cite{Lions-Sznitman-84} Theorem~3.1.
\begin{lem}
Assume that (G1) and (H1) hold true. Then for every $x \in \overline{G}$ there exists a unique adapted continuous couple of processes $\set{(X^x_t,K_t^x);t\geqslant 0}$ with values in $\overline{G} \times \mathbb{R}^+$ such that
\begin{equation}
\label{EDS reflechie}
\begin{array}{l}
X^x_t=x+\int_0^t b(X^x_s)ds+\int_0^t \sigma(X^x_s)dW_s+\int_0^t \nabla \phi(X^x_s) dK_s^x, \quad t \geqslant 0;\\
K_t^x = \int_0^t \mathbbm{1}_{X^x_s \in \partial G} dK_s^x, \quad K^x \textrm{ is non decreasing.}
\end{array}
\end{equation}
\end{lem}

This section is devoted to the following type of BSDEs with infinite horizon
\begin{equation}
\label{EDSRE}
Y^x_t=Y^x_T+ \int_t^T[\psi(X^x_s,Z^x_s)-\lambda]ds-\int_t^T Z^x_s dW_s, \quad 0 \leqslant t \leqslant T < +\infty,
\end{equation}
where $\lambda$ is a real number and is part of the unknowns of the problem and $\psi: \overline{G} \times \mathbb{R}^d \rightarrow \mathbb{R}$ verifies the following general assumptions:
\paragraph{(H2).} there exist $K_{\psi,x}\geqslant0$ and $K_{\psi,z}\geqslant0$ such that
\begin{equation*}
|\psi(x,z)-\psi(x^{\prime},z^{\prime})| \leqslant K_{\psi,x} |x-x^{\prime}|+K_{\psi,z}|z-z^{\prime}|, \quad \forall x,x^{\prime} \in \overline{G}, \, z,z^{\prime} \in \mathbb{R}^d.
\end{equation*}
We notice that $\psi(.,0)$ is continuous so there exists a constant $M_{\psi}$ verifying $\abs{\psi(.,0)} \leqslant M_{\psi}$. As in \cite{Fuhrman-Hu-Tessitore-07}, we start by considering an infinite horizon equation with strictly monotonic drift, namely, for $\alpha >0$, the equation
\begin{equation}
\label{EDSREapprochee}
Y_{t}^{x,\alpha}=Y_{T}^{x,\alpha}+ \int_t^T[\psi(X^x_s,Z_{s}^{x,\alpha})-\alpha Y_{s}^{x,\alpha}]ds-\int_t^T Z_{s}^{x,\alpha} dW_s, \quad 0 \leqslant t \leqslant T < +\infty.
\end{equation}
Existence and uniqueness have been first study by Briand and Hu in \cite{Briand-Hu-98} and then generalized by Royer in \cite{Royer-04}. They have established the following result:
\begin{lem}
\label{existence-unicite sol approchee}
Assume that (G1), (H1) and (H2) hold true. Then there exists a unique solution $(Y^{x,\alpha} , Z^{x,\alpha})$ to BSDE (\ref{EDSREapprochee}) such that $Y^{x,\alpha}$ is a bounded adapted continuous process and $Z^{x,\alpha} \in \mathcal{M}^2(\mathbb{R}^+,\mathbb{R}^d)$. Furthermore, $|Y_t^{x,\alpha}| \leqslant M_{\psi}/ \alpha$, $\mathbb{P}$-a.s. for all $t \geqslant 0$.
\end{lem}
We define 
$$v_{\alpha}(x):=Y_{0}^{x,\alpha}.$$ 
It is worth noting that $|v_{\alpha}(x)| \leqslant M_{\psi}/\alpha$ and uniqueness of solutions implies that $v_{\alpha}(X^x_t)=Y_{t}^{x,\alpha}$. The next step is to show that $v_{\alpha}$ is uniformly Lipschitz with respect to $\alpha$. Let 
\begin{equation*}
\eta:=\sup_{x,y \in \overline{G}, x\neq y} \left\{\dfrac{^t(x-y)(b(x)-b(y))}{|x-y|^{2}} +\dfrac{\textrm{Tr}[(\sigma(x)-\sigma(y))^t (\sigma(x)-\sigma(y))]}{2|x-y|^2}\right\}.
\end{equation*}
We will use the following assumption:
\paragraph{(H3).}
$\eta+K_{\psi,z}K_{\sigma}<0.$
\begin{rem}
When $\sigma$ is a constant function, (H3) becomes
$$\sup_{x,y \in \overline{G}, x\neq y} \left\{ \dfrac{^t(x-y)(b(x)-b(y))}{|x-y|^{2}}\right\}<0,$$
i.e. $b$ is dissipative.
\end{rem}

\begin{prop}
\label{valpha lipschitz}
Assume that (G1), (G2), (H1), (H2) and (H3) hold. Then we have, for all $\alpha>0$ and $x,x' \in \overline{G}$,
$$|v_{\alpha}(x)-v_{\alpha}(x^{\prime})|\leqslant \dfrac{K_{\psi,x}}{-\eta-K_{\psi,z}K_{\sigma}}|x-x^{\prime}|.$$
\end{prop}
\paragraph{Proof. }
We use a Girsanov argument due to P. Briand and Y. Hu in \cite{Briand-Hu-98}. Let $x,x^{\prime} \in \overline{G}$, we set $\tilde{Y}^{\alpha}:=Y^{x,\alpha}-Y^{x^{\prime},\alpha}$, $\tilde{Z}^{\alpha}:=Z^{x,\alpha}-Z^{x^{\prime},\alpha}$,
\begin{displaymath}
\label{beta}
 \beta(s)=\left\{
\begin{array}{cl}
 \dfrac{ \psi(X^{x^{\prime}}_s,Z_s^{x^{\prime},\alpha})-\psi(X^{x^{\prime}}_s,Z^{x,\alpha}_s) }{ |Z_s^{x^{\prime},\alpha}-Z_s^{x,\alpha}|^2 }\tr{(Z_s^{x^{\prime},\alpha}-Z_s^{x,\alpha})} & \textrm{if }Z_s^{x^{\prime},\alpha}-Z_s^{x,\alpha} \neq 0\\
0 &\textrm{otherwise,}\\
\end{array}
\right.
\end{displaymath}
\begin{eqnarray*} f_{\alpha}(s)&=&\psi(X^{x}_s,Z_s^{x,\alpha})-\psi(X^{x^{\prime}}_s,Z_s^{s,\alpha}),
\end{eqnarray*}
and $\tilde{W}_t=\int_0^t \beta_sds+W_t$. By hypothesis (H2), $\beta$ is a $\mathbb{R}^d$ valued adapted process bounded by $K_{\psi,z}$, so we are allowed to apply the Girsanov theorem: for all $T \in \mathbb{R}_+$ there exists a probability $\mathbb{Q}_T$ under which $(\tilde{W}_t)_{t \in [0,T]}$ is a Brownian motion. Then, from equation (\ref{EDSREapprochee}) we obtain
\begin{equation}
\label{EDSR approchée simplifiée}
\tilde{Y}_t^{\alpha}=\tilde{Y}_T^{\alpha}-\alpha\int_t^T\tilde{Y}_s^{\alpha}ds+\int_t^T f_{\alpha}(s) ds-\int_t^T \tilde{Z}_s^{\alpha} d\tilde{W}_s, \quad 0\leqslant t\leqslant T.
\end{equation}
Applying Itô's formula to $e^{-\alpha (s-t)}\tilde{Y}_s^{\alpha}$, we obtain
\begin{eqnarray*}
\tilde{Y}_t^{\alpha}&=&e^{-\alpha(T-t)}\tilde{Y}_T^{\alpha}+\int_t^T e^{-\alpha(s-t)}f_{\alpha}(s) ds-\int_t^T e^{-\alpha(s-t)}\tilde{Z}_s^{\alpha} d\tilde{W}_s\\
|\tilde{Y}_t^{\alpha}|&\leqslant&e^{-\alpha(T-t)} \mathbb{E}^{\mathbb{Q}_T}\Big[ |\tilde{Y}_T^{\alpha}|\Big|\mathcal{F}_t\Big] +\int_t^T e^{-\alpha(s-t)} \mathbb{E}^{\mathbb{Q}_T}\Big[ |f_{\alpha}(s)|\Big| \mathcal{F}_t\Big] ds\\
|\tilde{Y}_t^{\alpha}|&\leqslant&e^{-\alpha(T-t)} \mathbb{E}^{\mathbb{Q}_T}\Big[ |\tilde{Y}_T^{\alpha}|\Big|\mathcal{F}_t\Big]\\
&& +K_{\psi,x}\int_t^T e^{-\alpha(s-t)} \mathbb{E}^{\mathbb{Q}_T}\Big[ |X^x_s-X^{x^{\prime}}_s|^2\Big| \mathcal{F}_t\Big]^{1/2} ds.\\
\end{eqnarray*}
To conclude we are going to use the following lemma whose proof will be given after the proof of Theorem:
\begin{lem}
\label{X lipschitz}
Assume that (G1), (G2), (H1), (H2) and (H3) hold. For all $0\leqslant t \leqslant s \leqslant T$, $$\mathbb{E}^{\mathbb{Q}_T}\Big[ |X^x_s-X^{x^{\prime}}_s|^2\Big|\mathcal{F}_t\Big] \leqslant e^{2(\eta+K_{\psi,z}K_{\sigma})(s-t)} |X^x_t-X^{x^{\prime}}_t|^2.$$
Furthermore, if $\sigma$ is constant then, for all $0\leqslant t \leqslant s $, we have
$$|X^x_s-X^{x^{\prime}}_s| \leqslant e^{\eta(s-t)} |X^x_t-X^{x^{\prime}}_t|.$$
\end{lem}
From the last inequality, we deduce
\begin{eqnarray*}
|\tilde{Y}_t^{\alpha}|&\leqslant&e^{-\alpha(T-t)} \mathbb{E}^{\mathbb{Q}_T}\Big[ |\tilde{Y}_T^{\alpha}|\Big|\mathcal{F}_t\Big] +K_{\psi,x}|X^x_t-X^{x^{\prime}}_t|\int_t^T e^{(-\alpha+\eta+K_{\psi,z}K_{\sigma})(s-t)}ds,
\end{eqnarray*}
which implies
\begin{eqnarray*}
|\tilde{Y}_t^{\alpha}|&\leqslant&e^{-\alpha(T-t)}\dfrac{M_{\psi}}{\alpha} +K_{\psi,x}\dfrac{\big[1-e^{(-\alpha+\eta+K_{\psi,z}K_{\sigma})(T-t)}\big]}{\alpha-\eta-K_{\psi,z}K_{\sigma}}|X^x_t-X^{x^{\prime}}_t|.
\end{eqnarray*}
Finally, let $T \rightarrow +\infty$ and the claim follows by setting $t=0$. \cqfd
\paragraph{Proof of Lemma \ref{X lipschitz}. }
Let us apply Itô's formula to $e^{-2(\eta+K_{\psi,z}K_{\sigma})(s-t)}|X^x_s-X^{x^{\prime}}_s|^2$:
\begin{eqnarray*}
&&e^{-2(\eta+K_{\psi,z}K_{\sigma})(s-t)}|X^x_s-X^{x^{\prime}}_s|^2 = |X^x_t-X^{x^{\prime}}_t|^2\\
&&\qquad \qquad +2\int_t^s e^{-2(\eta+K_{\psi,z}K_{\sigma})(u-t)}\Big[ ^t(X^x_u-X^{x^{\prime}}_u)(b(X^x_u)-b(X^{x^{\prime}}_u))du\\
&&\qquad \qquad +\dfrac{1}{2}\textrm{Tr}[(\sigma(X^x_u)-\sigma(X^{x^{\prime}}_u))^t (\sigma(X^x_u)-\sigma(X^{x^{\prime}}_u))]du\\
&&\qquad \qquad +^t(X^x_u-X^{x^{\prime}}_u)\nabla\phi(X^{x}_u)dK_u^{x}-^t(X^x_u-X^{x^{\prime}}_u)\nabla\phi(X^{x^{\prime}}_u)dK_u^{x^{\prime}}\\
&&\qquad \qquad +^t(X^x_u-X^{x^{\prime}}_u)(\sigma(X^x_u)-\sigma(X^{x^{\prime}}_u))(d\tilde{W}_u-\beta_u du)\\
&&\qquad \qquad -(\eta+K_{\psi,z}K_{\sigma})|X^x_u-X^{x^{\prime}}_u|^2du\Big].
\end{eqnarray*}
$\overline{G}$ is a convex set, so $^t(x-y)\nabla\phi(x)\leqslant 0$ for all $(x,y) \in \partial G \times \overline{G}$. Furthermore $|\beta_s|\leqslant K_{{\psi},z}$ and $\sigma$ is $K_{\sigma}$-Lipschitz. By the definition of $\eta$ we obtain,
\begin{eqnarray*}
&&e^{2(-\eta-K_{\psi,z}K_{\sigma})(s-t)}|X^x_s-X^{x^{\prime}}_s|^2 \leqslant |X^x_t-X^{x^{\prime}}_t|^2\\ &&+2\int_t^s e^{-2(\eta+K_{\psi,z}K_{\sigma})(s-t)}\Big[ ^t(X^x_s-X^{x^{\prime}}_s)(\sigma(X^x_s)-\sigma(X^{x^{\prime}}_s))\Big]d\tilde{W}_s.
\end{eqnarray*}
Taking the conditional expectation of the inequality we get the first result. To conclude, the stochastic integral is a null function when $\sigma$ is a constant function. \cqfd
As in \cite{Fuhrman-Hu-Tessitore-07}, we now set $$\bar{v}_{\alpha}(x)=v_{\alpha}(x)-v_{\alpha}(0),$$
then we have $|\bar{v}_{\alpha}(x)|\leqslant \frac{K_{\psi,x}}{-\eta-K_{\psi,z}K_{\sigma}}|x|$ for all $x \in \overline{G}$ and all $\alpha >0$, according to Proposition~\ref{valpha lipschitz}. Moreover, $\alpha |v_{\alpha}(0)| \leqslant M_{\psi}$ by Lemma~\ref{existence-unicite sol approchee}. Thus we can construct by a diagonal procedure a sequence $(\alpha_n)_{n \in \mathbb{N}}\searrow 0$ such that, for all $x \in \overline{G} \cap \mathbb{Q}^d$, $\bar{v}_{\alpha_n}(x) \rightarrow\bar{v} (x)$ and $\alpha_nv_{\alpha_n}(0) \rightarrow \bar{\lambda}$. Furthermore, $\bar{v}_{\alpha}$ is a $\frac{K_{\psi,x}}{-\eta-K_{\psi,z}K_{\sigma}}$-Lipschitz function uniformly with respect to $\alpha$. So $\bar{v}$ can be extended to a $\frac{K_{\psi,x}}{-\eta-K_{\psi,z}K_{\sigma}}$-Lipschitz function defined on the whole  $\overline{G}$, thereby $\bar{v}_{\alpha_n}(x) \rightarrow\bar{v} (x)$ for all $x \in \overline{G}$. Thanks to this construction, we obtain the following theorem which can be proved in the same way as that of Theorem~4.4 in \cite{Fuhrman-Hu-Tessitore-07}.

\begin{thm}[Existence of a solution]
\label{existence convexe}
Assume that (G1), (G2), (H1), (H2) and (H3) hold. Let $\bar{\lambda}$ be the real number and $\bar{v}$ the function constructed previously. We define $\bar{Y}^x_t:=\bar{v}(X^x_t)$. Then, there exists a process $\bar{Z}^x \in \mathcal{M}^2(\mathbb{R}^+,\mathbb{R}^d)$ such that $\mathbb{P}-a.s.$ $(\bar{Y}^x,\bar{Z}^x,\bar{\lambda})$ is a solution of the EBSDE (\ref{EDSRE}) for all $x \in \overline{G}$. Moreover there exists a measurable function $\bar{\zeta}:\mathbb{R}^d \rightarrow \mathbb{R}$ such that $\bar{Z}^x_t=\bar{\zeta}(X^x_t)$.
\end{thm}

We remark that the solution to EBSDE (\ref{EDSRE}) is not unique. Indeed the equation is invariant with respect to addition of a constant to $Y$. However we have a result of uniqueness for $\lambda$.
\begin{thm}[Uniqueness of $\lambda$]
\label{unicite lambda g=0}
Assume that (G1), (H1) and (H2) hold. Let $(Y,Z,\lambda)$ a solution of EBSDE (\ref{EDSRE}). Then $\lambda$ is unique amongst solutions $(Y,Z,\lambda)$ such that $Y$ is a bounded continuous adapted process and $Z \in \mathcal{M}^2(\mathbb{R}^+,\mathbb{R}^d)$.
\end{thm}
\paragraph{Proof. }
We consider $(Y,Z,\lambda)$ and $(Y^{\prime},Z^{\prime},\lambda^{\prime})$ two solutions of the EBSDE (\ref{EDSRE}). Let $\tilde{\lambda}=\lambda^{\prime}-\lambda$, $\tilde{Y}=Y^{\prime}-Y$ and $\tilde{Z}=Z^{\prime}-Z$. We have, for all $T \in \mathbb{R}_+^*$,
$$\tilde{\lambda}=T^{-1}\left[ \tilde{Y}_T -\tilde{Y}_0\right]+T^{-1}\int_0^T \tilde{Z}_t\beta_tdt-T^{-1} \int_0^T \tilde{Z}_t dW_t$$
with
\begin{equation}
\label{beta unicite lambda}
 \beta_s=\left\{
\begin{array}{cl}
 \dfrac{ \psi(X^{x}_s,Z^{\prime}_s)-\psi(X^{x}_s,Z_s) }{ |Z^{\prime}_s-Z_s|^2 }\tr{(Z^{\prime}_s-Z_s)} & \textrm{if }Z^{\prime}_s-Z_s \neq 0\\
0 &\textrm{elsewhere.}\\
\end{array}
\right.
\end{equation}
$\beta$ is bounded: by the Girsanov theorem there exists a probability measure $\mathbb{Q}_T$ under which $(\tilde{W}_t=W_t-\int_0^t\beta_s ds)_{t\in [0,T]}$ is a Brownian motion. Computing the expectation with respect to $\mathbb{Q}_T$ we obtain
$$\tilde{\lambda}=T^{-1} \mathbb{E}^{\mathbb{Q}_T}\left[ \tilde{Y}_T -\tilde{Y}_0\right] \leqslant \frac{C}{T},$$
because $\tilde{Y}$ is bounded. So we can conclude the proof by letting $T\rightarrow +\infty$. \cqfd


To conclude this section we will show a proposition that will be usefull later.
\begin{prop}
 Assume that (G1), (H1) hold, $G$ is a bounded set and $\eta <0$. Then there exists a unique invariant measure $\nu$ for the process $(X_t)_{t\geqslant 0}$.
\end{prop}
\paragraph{Proof. }
The existence of an invariant measure $\nu$ for the process $(X_t)_{t\geqslant 0}$ is already stated in \cite{Skorokhod-89}, Theorem~1.21. Let $\nu$ and $\nu'$ two invariant measures and $X_0 \sim \nu$, $X'_0 \sim \nu'$ which are independent random variables of $(W_t)_{t\geqslant 0}$. For all $f \in \mathcal{C}_{lip}(\mathbb{R}^d)$ we have
$$|\mathbb{E}[f(X_0)]-\mathbb{E}[f(X'_0)]|=|\mathbb{E}[f(X^{X_0}_s)-f(X^{X'_0}_s)]|\leqslant K_f\mathbb{E}\Big[|X^{X_0}_s-X^{X'_0}_s|^2\Big]^{1/2},$$
with $K_f$ the Lipschitz constant of $f$. We are able to apply Lemma~\ref{X lipschitz} with $\psi=0$: for all $s \in \mathbb{R}^+$,
$$|\mathbb{E}[f(X_0)]-\mathbb{E}[f(X'_0)]|\leqslant K_fe^{-\eta s}\mathbb{E}\Big[|X_0-X'_0|^2\Big]^{1/2} \stackrel{s \rightarrow +\infty}{\longrightarrow} 0.$$
Then the claim ends by use of a density argument and the monotone class theorem. \cqfd

\section{EBSDEs with non-zero Neumann boundary conditions}
We are now interested in EBSDEs with non-zero Neumann boundary conditions: we are looking for solutions to the following type of BSDEs, for all  $0 \leqslant t \leqslant T < +\infty$,
\begin{equation}
\label{EDSRE g non nulle}
Y^x_t=Y^x_T+ \int_t^T[\psi(X^x_s,Z^x_s)-\lambda]ds+\int_t^T[g(X^x_s)-\mu]dK^x_s-\int_t^T Z^x_s dW_s,
\end{equation}
where $\lambda$ is a parameter, $\mu$ is part of the unknowns of the problem, $\psi$ still verifies (H2) and $g: \overline{G} \rightarrow \mathbb{R}$ verifies the following general assumption:
\paragraph{(F1). }
$g \in \mathcal{C}^2_{lip}(\overline{G})$.
\paragraph{}
Moreover we use extra assumption on $\phi$:
\paragraph{(G3). }
$\phi \in \mathcal{C}^2_{lip}(\mathbb{R}^d)$.
\paragraph{}
In this situation we will say that $(Y,Z,\mu)$ is a solution of EBSDE (\ref{EDSRE g non nulle}) with $\lambda$ fixed. But, due to our proof strategy, we will study firstly a modified problem where $\mu$ is a parameter and $\lambda$ is part of the unknowns. In this case, we will say that $(Y,Z,\lambda)$ is a solution of EBSDE (\ref{EDSRE g non nulle}) with $\mu$ fixed. We establish the following result of existence:
\begin{thm}[Existence of a solution]
\label{existence solution gdiff0}
Assume that (G1), (G2), (G3), (H1), (H2), (H3) and (F1) hold true. Then for any $\mu \in \mathbb{R}$ there exist $\lambda \in \mathbb{R}$, $v \in \mathcal{C}^0_{lip}(\overline{G})$, $\zeta:\mathbb{R}^d \rightarrow \mathbb{R}$ a measurable function such that, if we define $Y^x_t:=v(X^x_t)$ and $Z^x_t:=\zeta(X^x_t)$ then $Z^x \in \mathcal{M}^2(\mathbb{R}^+,\mathbb{R}^d)$ and $\mathbb{P}-a.s.$ $(Y^x,Z^x,\lambda)$ is a solution of EBSDE (\ref{EDSRE g non nulle}) with $\mu$ fixed, for all $x \in \overline{G}$.
\end{thm}
\paragraph{Proof. }
Our strategy is to modify EBSDE (\ref{EDSRE g non nulle}) in order to apply Theorem \ref{existence convexe}. According to the Theorem 3.2 of \cite{Ladyzenkaja-68} there exists $\alpha \in \mathbb{R}$ and $\tilde{v} \in \mathcal{C}^2_{lip}(\overline{G})$ such that 
\begin{displaymath}
\left\{
\begin{array}{ll}
\bigtriangleup \tilde{v}-\alpha \tilde{v}=0 & \forall x \in G\\
\frac{\partial \tilde{v}}{\partial n}(x)+g(x)=\mu, & \forall x \in \partial G.
\\
\end{array}
\right.
\end{displaymath}
We set $\tilde{Y}^x_t=\tilde{v}(X^x_t)$ and $\tilde{Z}^x_t=\tr{\nabla}\tilde{v}(X^x_t)\sigma(X^x_t)$. These processes verify for all $0 \leqslant t \leqslant T < +\infty$,
$$\tilde{Y}^x_t= \tilde{Y}^x_T-\int_t^T \mathcal{L}\tilde{v}(X^x_s)ds+\int_t^T[g(X^x_s)-\mu]dK_s^x-\int_t^T\tilde{Z}^x_sdW_s.$$
We now consider the following EBSDE with infinite horizon:
\begin{equation}
\label{EDSREmodifiee}
\bar{Y}^x_t=\bar{Y}^x_T+ \int_t^T[\bar{\psi}(X^x_s,\bar{Z}^x_s)-\lambda]ds-\int_t^T \bar{Z}^x_s dW_s, \quad 0 \leqslant t \leqslant T < +\infty,
\end{equation}
with $\bar{\psi}(x,z)=\mathcal{L}\tilde{v}(x)+\psi(x,z+\tr{\nabla}\tilde{v}(x)\sigma(x))$.
Since derivatives of $\tilde{v}$, $\sigma$ and $\psi$ are Lipschitz functions, there exists a constant $K_{\tilde{\psi},x}$ such that we have for all $x,x' \in \overline{G}$ and $z,z' \in \mathbb{R}^d$
$$|\tilde{\psi}(x,z)-\tilde{\psi}(x',z')|\leqslant K_{\tilde{\psi},x}|x-x'| + K_{\psi,z}|z-z'| \quad .$$
So we are able to apply Theorem \ref{existence convexe}: there exists $\bar{\lambda} \in \mathbb{R}$, $\bar{v} \in \mathcal{C}^0_{lip}(\overline{G})$ and $\bar{\xi}:\mathbb{R}^d \rightarrow \mathbb{R}$ a measurable function such that $(\bar{Y}^x:=\bar{v}(X^x),\bar{Z}^x:=\bar{\xi}(X^x),\bar{\lambda})$ is a solution of EBSDE (\ref{EDSREmodifiee}). We set
\begin{displaymath}
\begin{array}{l}
Y^x_t:=\tilde{Y}^x_t+\bar{Y}^x_t=\tilde{v}(X^x_t)+\bar{v}(X^x_t),\\
Z^x_t:=\tilde{Z}^x_t+\bar{Z}^x_t=\tr{\nabla}\tilde{v}(X^x_t)\sigma(X^x_t)+\bar{\xi}(X^x_t).
\\
\end{array}
\end{displaymath}
Then $(Y^x,Z^x,\bar{\lambda})$ is a solution of EBSDE (\ref{EDSRE g non nulle}) linked to $\mu$. \cqfd
We have also a result of uniqueness for $\lambda$ that can be shown exactly as Theorem \ref{unicite lambda g=0}:
\begin{thm}[Uniqueness of $\lambda$]
\label{unicite lambda gdiff0}
Assume that (G1), (H1) and (H2) hold. Let $(Y,Z,\lambda)$ a solution of EBSDE (\ref{EDSRE g non nulle}) with $\mu$ fixed. Then $\lambda$ is unique among solutions $(Y,Z,\lambda)$ such that $Y$ is a bounded continuous adapted process and $Z \in \mathcal{M}^2(\mathbb{R}^+,\mathbb{R}^d)$.
\end{thm}
Thanks to the uniqueness we can define the map $\mu \mapsto \lambda(\mu)$ and study its properties.

\begin{prop}
\label{lambda(mu) decroissante continue}
Assume that (G1), (G2), (G3), (H1), (H2), (H3) and (F1) hold true. Then $\lambda(\mu)$ is a decreasing continuous function on $\mathbb{R}$.
\end{prop}
\paragraph{Proof. }
Let $(Y^x,Z^x,\lambda)$ and $(\tilde{Y}^x,\tilde{Z}^x,\tilde{\lambda})$ two solutions of (\ref{EDSRE g non nulle}) linked to $\mu$ and $\tilde{\mu}$. We set $\bar{Y}^x:=\tilde{Y}^x-Y^x$ and $\bar{Z}^x:=\tilde{Z}^x-Z^x$. These processes verify for all $T \in \mathbb{R}_+$
\begin{equation}
\label{difference de deux EDSRE}
\bar{Y}_0^x=\bar{Y}_T^x+\int_0^T \big[\psi(X^x_s,\tilde{Z}^x_s)-\psi(X^x_s,Z^x_s)\big]ds +[\lambda-\tilde{\lambda}]T+[\mu-\tilde{\mu}]K_T^x-\int_0^T\bar{Z}^x_s dW_s.
\end{equation}

As usual, we set
\begin{displaymath}
 \beta_s=\left\{
\begin{array}{cl}
 \dfrac{ \psi(X^x_s,\tilde{Z}^x_s)-\psi(X^x_s,Z^x_s) }{ |\tilde{Z}^x_s-Z^x_s|^2 }\tr(\tilde{Z}^x_s-Z^x_s) & \textrm{if }\tilde{Z}^x_s-Z^x_s \neq 0\\
0 &\textrm{otherwise,}\\
\end{array}
\right.,
\end{displaymath}
and $\tilde{W}_t=-\int_0^t \beta_sds+W_t$. According to the Girsanov theorem there exists a probability $\mathbb{Q}_T$ under which $(\tilde{W}_t)_{t \in [0,T]}$ is a Brownian motion. Then we have
\begin{equation}
\label{etude lambda(mu) sous proba modifiee}
\bar{Y}_0^x=\underbrace{\mathbb{E}^{\mathbb{Q}_T}\Big[\bar{Y}_T^x\Big]}_{\leqslant M} +[\lambda-\tilde{\lambda}]T+[\mu-\tilde{\mu}]\mathbb{E}^{\mathbb{Q}_T}\underbrace{\Big[K_T^x\Big]}_{\geqslant 0}.
\end{equation}
If we suppose that $\mu \leqslant \tilde{\mu}$ and $\lambda < \tilde{\lambda}$ then 
$$\bar{Y}_0^x \leqslant [\lambda-\tilde{\lambda}]T+M \stackrel{n \rightarrow +\infty}{\longrightarrow}-\infty$$
this is a contradiction. So $\mu \leqslant \tilde{\mu} \Rightarrow \lambda \geqslant \tilde{\lambda}$. To show the continuity of $\lambda$ we assume that $\abs{\tilde{\mu}-\mu} \leqslant \eps$ with $\eps >0$. Then
$$\abs{\tilde{\lambda}-\lambda}=\frac{1}{T}\abs{\mathbb{E}^{\mathbb{Q}_T}\Big[ \bar{Y}_0^x-\bar{Y}_T^x+[\tilde{\mu}-\mu]K_T^x\Big]}\leqslant \frac{2M}{T}+\frac{\eps}{T} \mathbb{E}^{\mathbb{Q}_T}\Big[ K_T^x \Big].$$
Let us now prove a lemma about the bound on $\mathbb{E}^{\mathbb{Q}_T}\Big[ K_t^x \Big]$.
\begin{lem}
\label{croissance sous lineaire de K}
There exists a constant $C$ such that
 $$\mathbb{E}^{\mathbb{Q}_T}\Big[ K_t^x \Big]\leqslant C(1+t), \quad \forall T \in \mathbb{R}^+, \forall t \in [0,T], \forall x \in \overline{G}.$$
\end{lem}
\paragraph{Proof of the lemma. }
Applying Itô's formula to $\phi(X_t^x)$ we have for all $t \in \mathbb{R}^+$ and all $x \in \overline{G}$
\begin{equation}
\label{formule exacte pour K}
K_t^x=\phi(X_t^x)-\phi(x)-\int_0^t \mathcal{L} \phi(X_s^x)ds-\int_0^t\tr{\nabla\phi}(X_s^x)\sigma(X_s^x)dW_s.
\end{equation}
Then
\begin{eqnarray*}
\mathbb{E}^{\mathbb{Q}_T}\Big[ K_t^x \Big] &=& \mathbb{E}^{\mathbb{Q}_T}\Big[ \phi(X_t^x)-\phi(x)-\int_0^t \mathcal{L} \phi(X_s^x)ds-\int_0^t\tr{\nabla\phi}(X_s^x)\sigma(X_s^x)(\beta_s ds +d\tilde{W}_s) \Big]\\
& \leqslant & \mathbb{E}^{\mathbb{Q}_T}\Big[ \underbrace{\abs{\phi(X_t^x)}}_{\leqslant C/2}+\underbrace{\abs{\phi(x)}}_{\leqslant C/2}+\int_0^t \underbrace{\abs{\mathcal{L} \phi(X_s^x)}}_{\leqslant C/2}ds + \int_0^t\underbrace{\abs{\tr{\nabla\phi}(X_s^x)\sigma(X_s^x)\beta_s}}_{\leqslant C/2} ds \Big]\\
& \leqslant & C(1+t).
\end{eqnarray*}
\cqfd
Let us return back to the proof of Proposition~\ref{lambda(mu) decroissante continue}. By applying Lemma~\ref{croissance sous lineaire de K} we obtain
$$\abs{\tilde{\lambda}-\lambda}\leqslant \frac{2M}{T}+\frac{T+1}{T}C\eps \stackrel{T \rightarrow +\infty}{\longrightarrow} C\eps.$$
The proof is therefore completed. \cqfd

To prove our second theorem of existence we need to introduce a further assumption.
\paragraph{(F2). }
\begin{enumerate}
 \item $|\psi|$ is bounded by $M_{\psi}$;
 \item $\mathbb{E}[\mathcal{L}\phi(X_0)]<0$ if $X_0\sim \nu$ with $\nu$ the invariant measure for the process $(X_t)_{t\geqslant 0}$.
\end{enumerate}
\begin{thm}[existence of a solution]
\label{existence solution gdiff0 mufixe}
Assume that (G1), (G2), (G3), (H1), (H2), (H3), (F1) and (F2) hold true. Then for any $\lambda \in \mathbb{R}$ there exists $\mu \in \mathbb{R}$, $v \in \mathcal{C}^0_{lip}(\overline{G})$, $\zeta:\mathbb{R}^d \rightarrow \mathbb{R}$ a measurable function such that, if we define $Y^x_t:=v(X^x_t)$ and $Z^x_t:=\zeta(X^x_t)$ then $Z^x \in \mathcal{M}^2(\mathbb{R}^+,\mathbb{R}^d)$ and $\mathbb{P}-a.s.$ $(Y^x,Z^x,\mu)$ is a solution of EBSDE (\ref{EDSRE g non nulle}) with $\lambda$ fixed, for all $x \in \overline{G}$. Moreover we have 
$$\abs{\lambda(\mu)-\lambda(0)-\mu \mathbb{E}[\mathcal{L}\phi(X_0)]} \leqslant 2 M_{\psi}.$$
\end{thm}

\paragraph{Proof. }
Let $(Y,Z,\lambda(\mu))$ and $(\tilde{Y},\tilde{Z},\lambda(0))$ two solutions of equation (\ref{EDSRE g non nulle}) linked to $\mu$ and $0$ respectively. Let $X_0\sim \nu$ independent of $(W_t)_{t\geqslant 0}$. Then, from equation (\ref{difference de deux EDSRE}), we deduce for all $T \in \mathbb{R}^+$
$$\mathbb{E}\Big[\bar{Y}^{X_0}_0-\bar{Y}^{X_0}_T-[\lambda(\mu)-\lambda(0)]T-\mu K_T^{X_0}\Big] = \mathbb{E}\Big[\int_0^T \psi(X^{X_0}_s,\tilde{Z}^{X_0}_s)-\psi(X^{X_0}_s,Z^{X_0}_s)ds \Big],$$
from which we deduce that
$$\abs{\mathbb{E}\Big[\bar{Y}_0^{X_0}-\bar{Y}_T^{X_0}\Big]-[\lambda(\mu)-\lambda(0)]T-\mu \mathbb{E}\Big[ K_T^{X_0}\Big]} \leqslant 2M_{\psi}T .$$
By using equation (\ref{formule exacte pour K}) we have
\begin{eqnarray*}
\mathbb{E}\Big[ K_T^{X_0}\Big] &=& \mathbb{E}\Big[ \phi(X_T^{X_0})-\phi(X_0)-\int_0^T \mathcal{L} \phi(X_s^{X_0})ds \Big]\\
 &=& -\int_0^T  \mathbb{E}\Big[ \mathcal{L} \phi(X_s^{X_0}) \Big]ds\\
 &=& -\mathbb{E}\Big[ \mathcal{L} \phi(X_0) \Big] T.
\end{eqnarray*}
Combining the last two relations, we get
$$\abs{\frac{\mathbb{E}\Big[\bar{Y}_0^{X_0}-\bar{Y}_T^{X_0}\Big]}{T}-[\lambda(\mu)-\lambda(0)]+\mu \mathbb{E}\Big[ \mathcal{L} \phi(X_0)\Big]} \leqslant 2M_{\psi}.$$
Thus letting $T \rightarrow +\infty$ we conclude that
$$\abs{\lambda(\mu)-\lambda(0)-\mu \mathbb{E}[\mathcal{L}\phi(X_0)]} \leqslant 2M_{\psi}.$$
So, we obtain
$$\lambda(\mu) \stackrel{\mu \rightarrow +\infty}{\longrightarrow} - \infty \quad \textrm{and} \quad \lambda(\mu) \stackrel{\mu \rightarrow -\infty}{\longrightarrow} + \infty.$$
Finally the result is a direct consequence of the intermediate value theorem. \cqfd

The hypothesis $\mathbb{E}[\mathcal{L}\phi(X_0)]<0$ say that the boundary has to be visited recurrently. When $\sigma$ is non-singular on $\overline{G}$ we show that this hypothesis is always verified.
\begin{prop}
\label{F2 verifiee si sigma inversible}
Assume that (G1), (G2) and (H1) hold true. We assume also that $\sigma(x)$ is non-singular for all $x \in \overline{G}$. Then for the invariant measure $\nu$ of the process $(X_t)_{t\geqslant 0}$ we have $\mathbb{E}[\mathcal{L}\phi(X_0)]<0$ if $X_0\sim \nu$.
\end{prop}
\paragraph{Proof. }
Let us take a random variable $X_0 \sim \nu$ independent of $(W_t)_{t\geqslant 0}$. Then $\mathbb{E}\Big[ K_T^{X_0}\Big] = -\mathbb{E}\Big[ \mathcal{L} \phi(X_0) \Big] T$, which implies that $\mathbb{E}\Big[ \mathcal{L} \phi(X_0) \Big] \leqslant 0$. If $\mathbb{E}[\mathcal{L}\phi(X_0)]=0$, then $\mathbb{P}$-a.s. $K_t^{X_0}=0$, for all $t \in \mathbb{R}^+$. So the process $X^{X_0}$ is the solution of the stochastic differential equation
\begin{equation}
\label{EDS non reflechie}
X^{X_0}_t=X_0+\int_0^t \tilde{b}(X^{X_0}_s)ds+\int_0^t \tilde{\sigma}(X^{X_0}_s)dW_s, \quad t \geqslant 0,
\end{equation}
with $\tilde{b}$ and $\tilde{\sigma}$ defined on $\mathbb{R}^d$ by $\tilde{\sigma}(x)=\sigma(\textrm{proj}_{\overline{G}}(x))$ and $\tilde{b}(x)=b(\textrm{proj}_{\overline{G}}(x))$.
But according to \cite{Hasminskii-80} (Corollary~2 of Theorem~7.1), the solution of equation~(\ref{EDS non reflechie}) is a recurrent Markov process on $\mathbb{R}^d$. Thus this process is particularly unbounded: we have a contradiction. \cqfd
When $\sigma$ is singular on $\overline{G}$ then (F2) is not necessarily verified.
\paragraph{Examples. }
\begin{itemize}
 \item Let $\overline{G}=B(0,1)$, $\phi(x)=\frac{1-|x|^2}{2}$, $b(x)=-x$ and $\sigma(x)=\left(\begin{array}{lll}
x_1 &  & 0 \\ 
 & \ddots &  \\ 
0 &  & x_d
\end{array}\right)$ on $\overline{G}$. Then $\delta_0$ is an invariant measure and $\mathcal{L}(\phi)(0)=0$. If we set $d=1$, $\psi=0$ and $g=0$ then solutions of the differential equation~(\ref{intro EDP}) without boundary condition are $\set{A_i+B_ix^3-\frac{2}{3}\lambda \ln|x|,(A_i,B_i)\in \mathbb{R}^2} $ on $[-1,0[$ and $]0,1]$. Thereby bounded continuous solutions are $\set{A-\frac{\mu}{3}|x|^3, A \in \mathbb{R}}$ and $\lambda(\mu)=0$.
 \item Let $\overline{G}=B(0,1)$, $\phi(x)=\frac{1-|x|^2}{2}$, $b(x)=-x$ and $\sigma(x)=\left(\begin{array}{ll}
I_k & 0  \\ 
0 & 0_{d-k}
\end{array}\right)$ on $\overline{G}$.

$F_k:=\set{x \in \mathbb{R}^d / x_{k+1}=...=x_d=0}\simeq \mathbb{R}^k$ is a stationary subspace for solutions of equation~(\ref{EDS reflechie}). Let $\nu_k$ an invariant measure on $\mathbb{R}^k$ for $\tilde{\phi}(x)=\frac{1-|x|^2}{2}$, $\tilde{b}(x)=-x$ and $\tilde{\sigma}(x)=I_k$. According to Proposition~\ref{F2 verifiee si sigma inversible}, $\mathbb{E}^{\nu_k}[\tilde{\mathcal{L}}(\tilde{\phi})]<0$. Then $\nu:=\nu_k \otimes \delta_{0_{\mathbb{R}^{d-k}}}$ is an invariant measure for the initial problem and $\mathbb{E}^{\nu}[\mathcal{L}(\phi)]<0$.
\end{itemize}

\paragraph{}

Theorem~\ref{existence solution gdiff0 mufixe} is not totally satisfactory for two reasons: we have not a result on the uniqueness of $\mu$ and $\psi$ is usually not bounded in optimal ergodic control problems. So we introduce another result of existence with different hypothesis.
\paragraph{(F2'). }
$-\mathcal{L}\phi(x) > |\tr{\nabla \phi} \sigma|_{\infty,\overline{G}} K_{\psi,z}, \quad \forall x \in \overline{G}.$

\begin{thm}[Existence and uniqueness of a solution 2]
\label{existence solution gdiff0 mufixe 2}
Assume that (G1), (G2), (G3), (H1), (H2), (H3), (F1) and (F2') hold true. Then for any $\lambda \in \mathbb{R}$ there exists $\mu \in \mathbb{R}$, $v \in \mathcal{C}^0_{lip}(\overline{G})$, $\zeta:\mathbb{R}^d \rightarrow \mathbb{R}$ a measurable function such that, if we define $Y^x_t:=v(X^x_t)$ and $Z^x_t:=\zeta(X^x_t)$ then $Z^x \in \mathcal{M}^2(\mathbb{R}^+,\mathbb{R}^d)$ and $\mathbb{P}-a.s.$ $(Y^x,Z^x,\mu)$ is a solution of EBSDE (\ref{EDSRE g non nulle}) with $\lambda$ fixed, for all $x \in \overline{G}$. Moreover $\mu$ is unique among solutions $(Y,Z,\mu)$ with $\lambda$ fixed such that $Y$ is a bounded continuous adapted process and $Z \in \mathcal{M}^2(\mathbb{R}^+,\mathbb{R}^d)$.
\end{thm}

\paragraph{Proof. }
Let $(Y,Z,\lambda(\mu))$ and $(\tilde{Y},\tilde{Z},\lambda(\tilde{\mu}))$ two solutions of equation (\ref{EDSRE g non nulle}) linked to $\mu$ and $\tilde{\mu}$. As in the proof of Proposition~\ref{lambda(mu) decroissante continue} we set $\bar{Y}^x:=\tilde{Y}^x-Y^x$ and $\bar{Z}^x:=\tilde{Z}^x-Z^x$. From equation~\ref{etude lambda(mu) sous proba modifiee}, we have:
\begin{equation*}
(\mu-\tilde{\mu})\mathbb{E}^{\mathbb{Q}_T}\Big[\frac{K_T^{x}}{T}\Big]=\frac{1}{T}\Big(\bar{Y}_0^x-\mathbb{E}^{\mathbb{Q}_T}\big[\bar{Y}_T^x\big]\Big) -(\lambda(\mu)-\lambda(\tilde{\mu})).
\end{equation*}
$\bar{Y}^x$ is bounded, so $\mathbb{E}^{\mathbb{Q}_T}\big[ K_T^{x}/T\big]$ has a limit $l_{\mu,\tilde{\mu}}\geqslant 0$ when $T\rightarrow +\infty$ and $\mu \neq \mu'$ such that
\begin{equation}
\label{l_mu,mu'}
 (\lambda(\mu)-\lambda(\tilde{\mu}))+(\mu -\tilde{\mu})l_{\mu,\tilde{\mu}}=0.
\end{equation}

By use of equation (\ref{formule exacte pour K}) we have
\begin{eqnarray*}
\mathbb{E}^{\mathbb{Q}_T}\Big[ K_T^{x}\Big] & = & \mathbb{E}^{\mathbb{Q}_T}\Big[ \phi(X_T^{x})-\phi(x)-\int_0^T \mathcal{L} \phi(X_s^{x})ds -\int_0^T \tr{\nabla\phi}(X_s^x)\sigma(X_s^x)\beta_s ds\Big]\\
\mathbb{E}^{\mathbb{Q}_T}\Big[ \frac{K_T^{x}}{T}\Big] & \geqslant & -\frac{2|\phi|_{\infty}}{T}+\Big[-\sup_{x \in \overline{G}} \mathcal{L}\phi-|\nabla \phi \sigma|_{\infty,\overline{G}} K_{\psi,z} \Big].\\
\end{eqnarray*}
We set $c=-\sup_{x \in \overline{G}} \mathcal{L}\phi-|\nabla \phi \sigma|_{\infty,\overline{G}} K_{\psi,z}$. Since hypothesis (F2') holds true, we have $c>0$ and $l_{\mu,\tilde{\mu}}\geqslant c>0$ when $\mu \neq \mu'$.
Thus, thanks to equation~(\ref{l_mu,mu'}),
$$\lambda(\mu) \stackrel{\mu \rightarrow +\infty}{\longrightarrow} - \infty \quad \textrm{and} \quad \lambda(\mu) \stackrel{\mu \rightarrow -\infty}{\longrightarrow} + \infty.$$
Once again the existence result is a direct consequence of the intermediate value theorem. Moreover, if $\lambda(\mu)=\lambda(\tilde{\mu})$ then $\mu=\tilde{\mu}$. \cqfd
\begin{rem}
By applying Lemma~\ref{croissance sous lineaire de K} we show that $\mathbb{E}^{\mathbb{Q}_T}\big[ K_T^{x}/T\big]$ is bounded. So we have:
$$0<c \leqslant l_{\mu,\tilde{\mu}} \leqslant C, \quad \forall \mu \neq \tilde{\mu}.$$
\end{rem}
\begin{rem}
If we interest in the second example dealt in this section we see that (F2') hold true when $k/2-1>K_{\psi,z}$. 
\end{rem}

\section{Study of reflected kolmogorov processes case}
In this section, we assume that $(X_t)_{t\geqslant 0}$ is a reflected Kolmogorov process. The aim is to obtain an equivalent to Theorem~\ref{existence solution gdiff0 mufixe 2} with a less restrictive hypothesis than (F2'). We set $\sigma=\sqrt{2}I$ and $b=-\nabla U$ where $U : \mathbb{R}^d\rightarrow \mathbb{R}$ verify the following assumptions:
\paragraph{(H4). }
$U \in \mathcal{C}^2(\mathbb{R}^d)$, $\nabla U$ is a Lipschitz function on $\mathbb{R}^d$ and $\nabla^2 U \geqslant cI$ with $c>0$.
\paragraph{}
We notice that (H4) implies (H3) and (H1). Moreover, without loss of generality, we use an extra assumption on $\phi$:
\paragraph{(G4). }
$\nabla \phi$ is a Lipschitz function on $\mathbb{R}^d$.
\paragraph{}
To study the reflected process we will introduce the related penalized process:
\begin{equation*}
 X_t^{n,x}=x-\int_0^t \nabla U_n(X_s^{n,x})ds + \sqrt{2}B_t, \quad t\geqslant 0, \quad x \in \mathbb{R}^d, \quad n \in \mathbb{N},
\end{equation*}
with $U_n=U+n\textrm{d}^2(.,\overline{G})$. According to \cite{Gegout-Petit-Pardoux-96}, $\textrm{d}^2(.,\overline{G})$ is twice differentiable and $\nabla^2 \textrm{d}^2(.,\overline{G}) \geqslant 0$. So, we have $\nabla^2 U_n \geqslant c I$. Let $\mathcal{L}_n$ the transition semigroup generator of $(X_t^n)_{t\geqslant 0}$ with domain $\mathbb{D}_2(\mathcal{L}_n)$ on $L^2(\nu_n)$ and $\nu_n$ its invariant measure given by
\begin{equation*}
 \nu_n(dx)=\frac{1}{N_n}\exp(-U_n(x))dx, \textrm{ with } N_n=\int_{\mathbb{R}^d} \exp(-U_n(x))dx.
\end{equation*}
\begin{prop}
\label{convergence en loi de nu_n}
$\mathbb{E}^{\nu_n}[f] \stackrel{n \rightarrow +\infty}{\longrightarrow} \mathbb{E}^{\nu}[f]$ for all Lipschitz functions $f$. Particularly, $\nu_n$ converge weakly to $\nu$.
\end{prop}
The proof is given in the appendix. We obtain a simple corollary:
\begin{cor}
\label{loi nu}
 $\nu(dx)=\frac{1}{N}\exp(-U(x))1_{x \in \overline{G}} dx, \textrm{ with } N=\int_{\overline{G}} \exp(-U(x))dx.$
\end{cor}
We now introduce a different assumption that will replace (F2'):
\paragraph{(F2''). }
$\left( \dfrac{\delta}{\sqrt{2c}} +\sqrt{2}|\nabla \phi|_{\infty, \overline{G}}\right)K_{\psi,z} < -\mathbb{E}^{\nu}[\mathcal{L}\phi],$

with $\delta=\sup_{x \in \overline{G}} (\tr{\nabla U}(x)x) - \inf_{x \in \overline{G}} (\tr{\nabla U}(x)x)$.

\begin{thm}[Existence and uniqueness of a solution 3]
\label{existence solution gdiff0 mufixe 3}
Theorem~\ref{existence solution gdiff0 mufixe 2} remains true if we assume that (G1), (G2), (G3), (G4), (H2), (H4), (F1) and (F2'') hold.
\end{thm}
\paragraph{Proof. }
If we use notations of the previous section, it is sufficient to show that there exists a constant $C>0$ such that $\lim_{T \rightarrow +\infty} \mathbb{E}^{\mathbb{Q}_T}\Big[ \frac{K_T^{X_0}}{T}\Big]\geqslant C$ for all $\mu \neq \tilde{\mu}$, where $X_0 \sim \nu$ is independent of $(W_t)_{t\geqslant 0}$. We set $\eps$ and define $A_T$ such that
$$\varepsilon \in \left]  \dfrac{\delta}{\sqrt{2c}} K_{\psi,z}, -\mathbb{E}[\mathcal{L}\phi(X_0)]-\sqrt{2}|\nabla \phi|_{\infty, \overline{G}} K_{\psi,z} \right[,$$
$$A_T:=\set{-\frac{1}{T}\int_0^T \mathcal{L}\phi(X_s^{X_0})ds \leqslant -\mathbb{E}[\mathcal{L}\phi(X_0)]-\eps},$$
with $X_0 \sim \nu$ and $T>0$. $\varepsilon$ is well defined thanks to hypothesis (F2'').
\begin{eqnarray*}
\mathbb{E}^{\mathbb{Q}_T}\Big[ \frac{K_T^{X_0}}{T}\Big] & = & \mathbb{E}^{\mathbb{Q}_T}\Big[ \frac{\phi(X_T^{X_0})}{T}-\frac{\phi(X_0)}{T}-\frac{1}{T}\int_0^T \mathcal{L} \phi(X_s^{X_0})ds\\
& & -\frac{\sqrt{2}}{T}\int_0^T \tr{\nabla\phi}(X_s^{X_0})\beta_s ds\Big]\\
& \geqslant & -\frac{2|\phi|_{\infty}}{T}+\mathbb{E}^{\mathbb{Q}_T}\Big[  (\mathbb{E}[-\mathcal{L}\phi(X_0)]-\varepsilon)1_{^cA_T}-|\mathcal{L}\phi|_{\infty,\overline{G}}1_{A_T}\Big]\\
& & -\sqrt{2}|\nabla \phi|_{\infty, \overline{G}} K_{\psi,z}\\
& \geqslant & -\frac{2|\phi|_{\infty}}{T}+ (\mathbb{E}[-\mathcal{L}\phi(X_0)]-\varepsilon)(1-\mathbb{Q}_T(A_T))-|\mathcal{L}\phi|_{\infty,\overline{G}}\mathbb{Q}_T(A_T)\\
& & -\sqrt{2}|\nabla \phi|_{\infty, \overline{G}} K_{\psi,z}.\\
\end{eqnarray*}
By using Hölder's inequality with $p>1$ and $q>1$ such that $1/p+1/q=1$ we obtain
\begin{eqnarray*}
\mathbb{Q}_T(A_T) &=& \mathbb{E}\left[ \exp\left(\int_0^T \beta_s dW_s-\frac{1}{2}\int_0^T |\beta_s|^2ds\right) 1_{A_T} \right]\\
&\leqslant& \mathbb{E}\left[ \exp \left( p\int_0^T \beta_s dW_s-\frac{p^2}{2}\int_0^T |\beta_s|^2ds+\frac{p(p-1)}{2}\int_0^T |\beta_s|^2ds \right) \right]^{1/p} \mathbb{P}(A_T)^{1/q}\\
&\leqslant& \exp \left( \frac{(p-1)}{2}K_{\psi,z}^2 T \right) \mathbb{P}(A_T)^{1-1/p}.\\
\end{eqnarray*}
To conclude we are going to use the following proposition which will be proved in the appendix thanks to Theorem~3.1 of \cite{Guillin-Leonard-Wu-Yao-07}:
\begin{prop}
\label{generalisation ineg grandes dev}
 Assume that (G1), (G2), (G3), (G4), (H1) and (H4) hold. Then
$$\mathbb{P}(A_T) \leqslant \exp \left( -\frac{c\varepsilon^2 T}{\delta^2}  \right).$$
\end{prop}
So
$$\mathbb{Q}_T(A_T) \leqslant \exp \left[ \underbrace{\left(\frac{p(p-1)}{2}K_{\psi,z}^2 -\frac{(p-1)c\varepsilon^2}{\delta^2} \right)}_{B_p} \frac{T}{p} \right].  $$
$B_p$ is a trinomial in $p$ that has two different real roots $1$ and $\frac{2c\varepsilon^2}{ \delta^2K_{\psi,z}^2}>1$ because $\varepsilon > \delta K_{\psi,z}/\sqrt{2c}$ by hypothesis (F2''). So we are able to find $p>1$ such that $B_p<0$. Then $\mathbb{Q}_T(A_T) \stackrel{T \rightarrow +\infty}{\longrightarrow} 0$ and 
$$\lim_{T \rightarrow +\infty} \mathbb{E}^{\mathbb{Q}_T}\Big[ \frac{K_T^{X_0}}{T}\Big] \geqslant -\mathbb{E}[\mathcal{L}\phi(X_0)]- 
\sqrt{2}|\nabla \phi|_{\infty, \overline{G}} K_{\psi,z}-\eps>0.$$\cqfd

\begin{rem}
All these results stay true if $\sigma(x)=\sqrt{2}\left(\begin{array}{ll}
I_k & 0  \\ 
0 & 0_{d-k}
\end{array}\right)$ and $F_k$, defined in the previous example, is a stationary subspace of $\nabla U$. We can even replace (F2'') by
$$\left( \sqrt{\frac{1}{2c}}\delta +\sqrt{2}|\nabla \phi|_{\infty, \overline{G}\cap F_k}\right)K_{\psi,z} < -\mathbb{E}^{\nu}[\mathcal{L}\phi],$$
with $\delta=\sup_{x \in \overline{G}\cap F_k} (\tr{\nabla U}(x)x) - \inf_{x \in \overline{G}\cap F_k} (\tr{\nabla U}(x)x)$. Indeed, as we see in the previous example, $\nu$ is nonzero at most on the set $\overline{G}\cap F_k$. So it is possible to restrict the process to the subspace $F_k$.
\end{rem}


\section[Prob. interpretation of the solution of an elliptic PDE]{Probabilistic interpretation of the solution of an elliptic PDE with linear Neumann boundary condition}
\label{interpretation probabiliste de la solution d'une EDP}
Consider the semi-linear elliptic PDE:
\begin{equation}
\label{EDP}
\left\{
\begin{array}{l}
\mathcal{L}v(x)+\psi(x,^t\nabla v(x)\sigma(x))=\lambda, \quad x \in G\\
\frac{\partial v}{\partial n}(x)+g(x)=\mu, \quad x \in \partial G,
\end{array}
\right.
\end{equation}
with
\begin{equation*}
\mathcal{L}f(x)=\dfrac{1}{2} \textrm{Tr}(\sigma(x)^t\sigma(x)\nabla^2f(x))+^tb(x)\nabla f(x).
\end{equation*}
We will prove now that $v$, defined in Theorem~\ref{existence solution gdiff0} or in Theorem~\ref{existence solution gdiff0 mufixe}, is a viscosity solution of PDE (\ref{EDP}). See e.g. \cite{Pardoux-Zhang-98} Definition~5.2 for the definition of a viscosity solution.
\begin{thm}
$v \in \mathcal{C}_{lip}^0(\overline{G})$, defined in Theorem~\ref{existence solution gdiff0} or in Theorem~\ref{existence solution gdiff0 mufixe}, is a viscosity solution of the elliptic PDE~(\ref{EDP}).
\end{thm}
\paragraph{Proof . }
It is a very standard proof that we can adapt easily from~\cite{Pardoux-Zhang-98}, Theorem 4.3. \cqfd
\begin{rem}
 With other hypothesis, uniqueness of solution $v$ is given by Barles and Da Lio in Theorem~4.4 of~\cite{Barles-DaLio-05}.
\end{rem}

If $\sigma$ is non-singular on $\overline{G}$ we notice that it is possible to jointly modify $b$ and $\psi$ without modify the PDE \ref{EDP}. We set $\tilde{b}(x)=b(x)-\xi x$ and $\tilde{\psi}(x,z)=\psi(x,z)+\xi z \sigma^{-1}(x) x$ for $\xi \in \mathbb{R}^+$. Then we are able to find a new hypothesis substituting (H3). We note $\tilde{\eta}$ the scalar $\eta$ corresponding to $\tilde{b}$.
\begin{prop}
 If $\eta+K_{\psi,z}K_{\sigma}<0$ or $K_{\sigma}\sup_{x \in \overline{G}} |\sigma^{-1}(x)x|<1$ then there exists $\xi\geqslant 0$ such that $\tilde{\eta}+K_{\tilde{\psi},z}K_{\sigma}<0$. In particular it is true when $\sigma$ is a constant function.
\end{prop}
\paragraph{Proof: }
It suffices to notice that $\tilde{\eta}=\eta-\xi$ and $K_{\tilde{\psi},z}\leqslant K_{\psi,z}+\xi\sup_{x \in \overline{G}} |\sigma^{-1}(x)x|$. So
\begin{equation*}
\tilde{\eta}+K_{\tilde{\psi},z}K_{\sigma}\leqslant\eta+K_{\psi,z}K_{\sigma}+\xi(K_{\sigma}\sup_{x \in \overline{G}} |\sigma^{-1}(x)x|-1).
\end{equation*}
\cqfd

\section{Optimal ergodic control}
Let $U$ be a separable metric space. We define a control $\rho$ as an $(\mathcal{F}_t)$-progressively measurable $U$-valued process. We introduce $R: U \rightarrow \mathbb{R}^d$ and $L: \mathbb{R}^d\times \mathbb{R}^{1\times d} \rightarrow \mathbb{R}$ two continuous functions such that, for some constants $M_R>0$ and $M_L>0$,
\begin{equation}
\label{hypotheses R et L}
|R(u)|\leqslant M_R, \quad |L(x,u)|\leqslant M_L, \quad |L(x,u)-L(x',u)|\leqslant c|x-x'|, \quad \forall u \in U,\,x,x' \in \mathbb{R}^d.
\end{equation}
Given an arbitrary control $\rho$ and $T>0$, we introduce the Girsanov density
$$\Gamma_T^{\rho}=\exp\left( \int_0^T R(\rho_s)dW_s-\frac{1}{2}\int_0^T |R(\rho_s)|^2ds \right)$$
and the probability $\mathbb{P}_T^{\rho}=\Gamma_T^{\rho}\mathbb{P}$ on $\mathcal{F}_T$.
Ergodic costs corresponding to a given control $\rho$ and a starting point $x \in \mathbb{R}^d$ are defined in the following way:
\begin{equation}
 \label{cout I}
I(x,\rho)=\limsup_{T \rightarrow +\infty} \frac{1}{T} \mathbb{E}^{\rho,T} \left[ \int_0^T L(X_s^x,\rho_s)ds+\int_0^T[g(X_s^x)-\mu]dK_s^x \right],
\end{equation}

\begin{equation}
 \label{cout J}
J(x,\rho)=\limsup_{T \rightarrow +\infty} \frac{1}{\mathbb{E}^{\rho,T}[K_T^x]} \mathbb{E}^{\rho,T} \left[ \int_0^T [L(X_s^x,\rho_s)-\lambda]ds+\int_0^Tg(X_s^x)dK_s^x \right] \mathbbm{1}_{\mathbb{E}^{\rho,T}[K_T^x]>0},
\end{equation}
where $\mathbb{E}^{\rho,T}$ denotes expectation with respect to $\mathbb{P}_T^{\rho}$. We notice that $W_t^{\rho}=W_t+\int_0^tR(\rho_s)ds$ is a Wiener process on  $[0,T]$ under $\mathbb{P}_T^{\rho}$. 


Our purpose is to minimize costs $I$ and $J$ over all controls. So we first define the Hamiltonian in the usual way
\begin{equation}
\label{psi comme hamiltonien}
 \psi(x,z)=\inf_{u \in U}\set{L(x,u)+zR(u)}, \quad x \in \mathbb{R}^d, z \in  \mathbb{R}^{1\times d},
\end{equation}
and we remark that if, for all $x,z,$ the infimum is attained in (\ref{psi comme hamiltonien}) then, according to Theorem~4 of \cite{MacShane-Warfield-67}, there exists a measurable function $\gamma: \mathbb{R}^d \times \mathbb{R}^{1\times d} \rightarrow U$ such that 
$$\psi(x,z)=L(x,\gamma(x,z))+zR(\gamma(x,z)).$$
We notice that $\psi$ is a Lipschitz function: hypothesis (H2) is verified with $K_{\psi,z}=M_R$.
\begin{thm}
\label{controle optimal I}
 Assume that hypothesis of Theorem~\ref{existence solution gdiff0} hold true. Let $(Y,Z,\lambda)$ a solution of (\ref{EDSRE g non nulle}) with $\mu$ fixed. Then the following holds:
\begin{enumerate}
 \item For arbitrary control $\rho$ we have $I(x,\rho)\geqslant\lambda$ and the equality holds if and only if $L(X_t^x,\rho_t)+Z_t^xR(\rho_t)=\psi(X_t^x,Z_t^x)$, $\mathbb{P}$-a.s. for almost every $t$.
 \item If the minimum is attained in (\ref{psi comme hamiltonien}) then the control $\overline{\rho}_t=\gamma(X_t^x,Z_t)$ verifies $I(x,\overline{\rho})=\lambda$.
\end{enumerate}
\end{thm}
\paragraph{Proof. }
This theorem can be proved in the same manner as that of Theorem~7.1 in \cite{Fuhrman-Hu-Tessitore-07} and we omit it.
\cqfd
\begin{rem}
\label{remarques sur cout optimal I}
 \begin{enumerate}
  \item If the minimum is attained in (\ref{psi comme hamiltonien}) then there exists an optimal feedback control given by the function $x \mapsto \gamma(x,\xi(x))$  where $(Y,\xi(X),\lambda)$ is the solution constructed in Theorem~\ref{existence solution gdiff0}.
  \item If limsup is changed into liminf in the definition (\ref{cout I}) of the cost, then the same conclusion hold, with the obvious modifications, and the optimal value is given by $\lambda$ in both cases.
 \end{enumerate}
\end{rem}

\begin{thm}
\label{controle optimal J}
 Assume that hypothesis of Theorem~\ref{existence solution gdiff0 mufixe 2} or Theorem~\ref{existence solution gdiff0 mufixe 3} hold true. 
Let $(Y,Z,\mu)$ a solution of (\ref{EDSRE g non nulle}) with $\lambda$ fixed. Then the following holds:
\begin{enumerate}
 \item For arbitrary control $\rho$ we have $J(x,\rho)\geqslant\mu$ and the equality holds if and only if $L(X_t^x,\rho_t)+Z_t^xR(\rho_t)=\psi(X_t^x,Z_t^x)$, $\mathbb{P}$-a.s. for almost every $t$.
 \item If the minimum is attained in (\ref{psi comme hamiltonien}) then the control $\overline{\rho}_t=\gamma(X_t^x,Z_t)$ verifies $J(x,\overline{\rho})=\mu$.
\end{enumerate}
\end{thm}
\paragraph{Proof. }
As $(Y,Z,\mu)$ is a solution of the EBSDE with $\lambda$ fixed, we have
\begin{eqnarray*}
 -dY_t^x &=& [\psi(X_t^x,Z_t^x)-\lambda]dt +[g(X_t^x)-\mu]dK_t^x - Z_t^x dW_t\\
  &=& [\psi(X_t^x,Z_t^x)-\lambda]dt +[g(X_t^x)-\mu]dK_t^x - Z_t^x dW_t^{\rho} -Z_t^xR(\rho_t)dt,\\
\end{eqnarray*}
from which we deduce that
\begin{eqnarray*}
 \mu\mathbb{E}^{\rho,T}[K_T^x] &=& \mathbb{E}^{\rho,T}\left[Y_T^x-Y_0^x\right] + \mathbb{E}^{\rho,T}\left[ \int_0^T [\psi(X_t^x,Z_t^x) - Z_t^xR(\rho_t) - L(X_t^x,\rho_t)] dt \right]\\
& & +\mathbb{E}^{\rho,T}\left[ \int_0^T [L(X_t^x,\rho_t)-\lambda] dt \right] + \mathbb{E}^{\rho,T} \left[\int_0^T g(X_t^x) dK_t^x \right].
\end{eqnarray*}
Thus
$$\mu\mathbb{E}^{\rho,T}[K_T^x] + \mathbb{E}^{\rho,T}\left[Y_0^x-Y_T^x\right] \leqslant \mathbb{E}^{\rho,T}\left[ \int_0^T [L(X_t^x,\rho_t)-\lambda] dt + \int_0^T g(X_t^x) dK_t^x \right].$$
To conclude we are going to use the following lemma that we will prove immediately after the proof of this theorem:
\begin{lem}
\label{E rho,T [KT] tend vers l'infini}
 Assume that hypothesis of Theorem~\ref{existence solution gdiff0 mufixe 2} or Theorem~\ref{existence solution gdiff0 mufixe 3} hold true. Then for all $x \in \overline{G}$
$$\lim_{T \rightarrow +\infty} \mathbb{E}^{\rho,T}[K_T^x]=+\infty.$$
\end{lem}
So, for $T>T_0$, $\mathbb{E}^{\rho,T}[K_T^x]>0$ and
$$\mu + \frac{\mathbb{E}^{\rho,T}\left[Y_0^x-Y_T^x\right]}{\mathbb{E}^{\rho,T}[K_T^x]} \leqslant \frac{1}{\mathbb{E}^{\rho,T}[K_T^x]}\mathbb{E}^{\rho,T}\left[ \int_0^T [L(X_t^x,\rho_t)-\lambda] dt + \int_0^T g(X_t^x) dK_t^x \right].$$
Since $Y$ is bounded we finally obtain
$$\mu \leqslant \limsup_{T \rightarrow +\infty} \frac{1}{\mathbb{E}^{\rho,T}[K_T^x]}\mathbb{E}^{\rho,T}\left[ \int_0^T [L(X_t^x,\rho_t)-\lambda] dt + \int_0^T g(X_t^x) dK_t^x \right]=J(x,\rho).$$
Similarly, if $L(X_t^x,\rho_t)+Z_t^xR(\rho_t)=\psi(X_t^x,Z_t^x),$
$$\mu\mathbb{E}^{\rho,T}[K_T^x] + \mathbb{E}^{\rho,T}\left[Y_0^x-Y_T^x\right] = \mathbb{E}^{\rho,T}\left[ \int_0^T [L(X_t^x,\rho_t)-\lambda] dt + \int_0^T g(X_t^x) dK_t^x \right],$$
and the claim holds.
\cqfd
\paragraph{Proof of Lemma~\ref{E rho,T [KT] tend vers l'infini}. }
Firstly we assume that hypothesis of Theorem~\ref{existence solution gdiff0 mufixe 2} hold true. As in the proof of this theorem, we have by using equation~(\ref{formule exacte pour K}),
\begin{eqnarray*}
\mathbb{E}^{\rho,T}\Big[ K_T^{x}\Big] & = & \mathbb{E}^{\rho,T}\left[ \phi(X_T^{x})-\phi(x)-\int_0^T \mathcal{L} \phi(X_s^{x})ds -\int_0^T \tr{\nabla\phi}(X_s^x)\sigma(X_s^x)R(\rho_s) ds\right],
\end{eqnarray*}
from which we deduce that
\begin{eqnarray*}
\mathbb{E}^{\rho,T}\left[ \frac{K_T^{x}}{T}\right] & \geqslant & -\frac{2|\phi|_{\infty}}{T}+\Big[-\sup_{x \in \overline{G}} \mathcal{L}\phi(x)-|\nabla \phi \sigma|_{\infty,\overline{G}} M_R \Big].\\
\end{eqnarray*}
Thanks to hypothesis (F2') we have
\begin{eqnarray*}
\mathbb{E}^{\rho,T}\left[ \frac{K_T^{x}}{T}\right] & \geqslant & \frac{1}{2}\Big[-\sup_{x \in \overline{G}} \mathcal{L}\phi(x)-|\nabla \phi \sigma|_{\infty,\overline{G}} M_R \Big]>0, \qquad \forall T>T_0,
\end{eqnarray*}
and the claim is proved. We now assume that hypothesis of Theorem~\ref{existence solution gdiff0 mufixe 3} hold true. Let $X_0 \sim \nu$ be a random variable independent of $(W_t)_{t\geqslant 0}$ and $\nu$ the invariant measure of $(X_t)_{t\geqslant 0}$. Exactly as in the proof of Theorem~\ref{existence solution gdiff0 mufixe 3} we are able to show that $\mathbb{E}^{\rho,T}\left[ K_T^{X_0}/T\right] \geqslant C >0$ for all $T>T_0$ by replacing $\beta$ with $R(\rho)$. On the other hand, for all $x \in \overline{G}$ and $T \in \mathbb{R}_+^*$, we have
\begin{eqnarray*}
 \abs{\frac{\mathbb{E}^{\rho,T}\Big[ K_T^{X_0}\Big]-\mathbb{E}^{\rho,T}\Big[ K_T^{x}\Big]}{T}} &\leqslant & \frac{4|\phi|_{\infty}}{T}+\frac{1}{T}\mathbb{E}^{\rho,T} \int_0^T|\mathcal{L}\phi(X_s^{X_0})-\mathcal{L}\phi(X_s^{x})|ds\\
& & +\frac{1}{T}\mathbb{E}^{\rho,T} \int_0^T|\tr{\nabla\phi}(X_s^{X_0})\sigma(X_s^{X_0})-\tr{\nabla\phi}(X_s^{x})\sigma(X_s^{x})||R(\rho_s)|ds\\
\end{eqnarray*}
Since $\mathcal{L}\phi$ and $\tr{\nabla\phi}\sigma$ are Lipschitz functions, we obtain
\begin{eqnarray*}
\abs{\frac{\mathbb{E}^{\rho,T}\Big[ K_T^{X_0}\Big]-\mathbb{E}^{\rho,T}\Big[ K_T^{x}\Big]}{T}} &\leqslant & \frac{4|\phi|_{\infty}}{T}+\frac{K_{\mathcal{L}\phi}}{T}\mathbb{E}^{\rho,T} \int_0^T|X_s^{X_0}-X_s^{x}|ds\\
& & +\frac{M_RK_{\tr{\nabla\phi}\sigma}}{T}\mathbb{E}^{\rho,T} \int_0^T|X_s^{X_0}-X_s^{x}|ds.\\
\end{eqnarray*}
Exactly as in Lemma~\ref{X lipschitz} we are able to show that for all $s\geqslant 0$
$$\mathbb{E}^{\rho,T}\left[|X_s^{X_0}-X_s^{x}|^2\right] \leqslant e^{2(\eta+M_R K_{\sigma})s} \mathbb{E}^{\rho,T}\left[|X_0-x|^2\right].$$
Finally,
\begin{eqnarray*}
\abs{\frac{\mathbb{E}^{\rho,T}\Big[ K_T^{X_0}\Big]-\mathbb{E}^{\rho,T}\Big[ K_T^{x}\Big]}{T}} &\leqslant & \frac{K_{\mathcal{L}\phi}+M_RK_{\tr{\nabla\phi}\sigma}}{T} \mathbb{E}^{\rho,T}\left[|X_0-x|^2\right]^{1/2} \int_0^T e^{(\eta+M_RK_{\sigma})s} ds\\
& & + \frac{4|\phi|_{\infty}}{T}\\
 &\leqslant & \frac{K_{\mathcal{L}\phi}+M_RK_{\tr{\nabla\phi}\sigma}}{T} \mathbb{E}^{\rho,T}\left[|X_0-x|^2\right]^{1/2} \frac{1-e^{(\eta+M_RK_{\sigma})T}}{-\eta-M_RK_{\sigma}}\\
& & +\frac{4|\phi|_{\infty}}{T}.
\end{eqnarray*}
Since hypothesis (H3) holds true, $\eta+M_RK_{\sigma} <0$ and so 
$$\lim_{T \rightarrow +\infty} \abs{\frac{\mathbb{E}^{\rho,T}\Big[ K_T^{X_0}\Big]-\mathbb{E}^{\rho,T}\Big[ K_T^{x}\Big]}{T}}=0.$$
Thus, for all $x \in \overline{G}$ there exists $T_0\geqslant 0$ such that
$$\mathbb{E}^{\rho,T}\left[K_T^{x}/T\right] \geqslant \frac{1}{2}\mathbb{E}^{\rho,T}\left[K_T^{X_0}/T\right] \geqslant c/2 >0$$
and the claim follows.
\cqfd

\begin{rem}
 Remarks~\ref{remarques sur cout optimal I} remains true for Theorem~\ref{controle optimal J}.
\end{rem}

\section[EBSDEs on a non-convex bounded set]{Some additional results: EBSDEs on a non-convex bounded set}
\label{section G non convexe}
In previous sections we have supposed that $G$ was a bounded convex set. We shall substitute hypothesis (G2) by this one:
\paragraph{(G2'). } $G$ is a bounded subset of $\mathbb{R}^d$.
\paragraph{}
In this section we suppose also that $\sigma$ is a constant function. At last, we set
$$\alpha=\sup_{x \in co(\bar{G})} \sup_{|y|=1} ( ^t y \nabla^2 \phi(x) y)$$
with $co(\bar{G})$ the convex hull of $\bar{G}$. Without loss of generality we assume that $\alpha>0$. Indeed, $\alpha \leqslant 0$ if and only if $\phi$ is concave which implies $\bar{G}$ is a convex set. In previous sections hypothesis (G2) has been used to prove Lemma~\ref{X lipschitz} so we will modify it:
\begin{lem}
Assume (G1), (G2'), (H1), (H2) hold true and $\sigma$ is a constant function. Let
\begin{eqnarray*}
\theta:=& &\sup_{x,y \in \bar{G},x\neq y, z,z' \in \R^d, z \neq z'} \bigg\{ 2\dfrac{\tr(x-y)(b(x)-b(y))}{|x-y|^{2}}\\
& &-\alpha \tr(\nabla \phi(x) + \nabla \phi(y))\sigma \beta(x,y,z,z')\\
& & -\frac{\alpha}{2} \trace \left(\nabla^2\phi(x)\sigma\tr \sigma+\nabla^2\phi(y)\sigma\tr\sigma\right)-\alpha\tr\nabla\phi(x)b(x)-\alpha\tr\nabla\phi(y)b(y)\\
& &+\alpha^2 \Big(\tr\nabla\phi(x)+\tr\nabla\phi(y)\Big)\sigma\tr\sigma\Big(\nabla\phi(x)+\nabla\phi(y)\Big) \bigg\},
\end{eqnarray*}
with $(z-z')\beta(x,y,z,z')=\big(\psi(x,z)+\psi(y,z)-\psi(x,z')-\psi(y,z')\big)/2$.
Then there exists a constant $M$ which depends only on $\phi$ and such that for all $0\leqslant t \leqslant s \leqslant n$,
$$\mathbb{E}^{\mathbb{Q}_n}\Big[ |X^x_s-X^{x^{\prime}}_s|^2\Big|\mathcal{F}_t\Big] \leqslant Me^{\theta(s-t)} |X^x_t-X^{x^{\prime}}_t|^2.$$
\end{lem}
\begin{rem}
$\beta$ exists, we can take 
\begin{displaymath}
 \beta=\left\{
\begin{array}{cl}
 \dfrac{ \psi(x,z')+\psi(y,z')-\psi(y,z)-\psi(x,z) }{ 2|z'-z|^2 }\tr{(z'-z)} & \textrm{if }z \neq z' \\
0 &\textrm{otherwise,}\\
\end{array}
\right.
\end{displaymath}
but there is not uniqueness. We have $|\beta|\leqslant K_{\psi,z}$ yet.
\end{rem}

\paragraph{Proof. }
Firstly we show an elementary lemma.
\begin{lem}
\label{lemme non convexe}
$\forall x \in \bar{G}$, $\forall y \in \partial G$ we have
$$-\alpha|x-y|^2+2\tr(y-x)\nabla\phi(y) \leqslant 0.$$
\end{lem}
\paragraph{Proof. }
Let $x\in \bar{G}$ and $y \in \partial G$. According to Taylor-Lagrange theorem there exists $t\in ]0,1[$ such that
$$\phi(x)=\phi(y)+\tr(x-y)\nabla\phi(y)+\frac{1}{2}\tr(x-y)\nabla^2\phi(tx+(1-t)(y-x))(x-y).$$
$\phi(x)\geqslant 0$, $\phi(y)=0$ and the claim easily follows. \cqfd
As in Lions and Sznitman~\cite{Lions-Sznitman-84} page 524,  using Itô's formula, we develop the semimartingale
$e^{-\theta u}e^{-\alpha(\phi(X^x_u)+\phi(X^{x'}_u))} |X^x_u-X^{x'}_u|^2,$
which leads us to
\begin{displaymath}
\begin{array}{l}
d\Big(e^{-\theta u}e^{-\alpha(\phi(X^x_u)+\phi(X^{x'}_u))} |X^x_u-X^{x'}_u|^2\Big)=\\
\qquad \qquad -\theta e^{-\theta u}e^{-\alpha(\phi(X^x_u)+\phi(X^{x'}_u))} |X^x_u-X^{x'}_u|^2du\\
\qquad \qquad +2e^{-\theta u}e^{-\alpha(\phi(X^x_u)+\phi(X^{x'}_u))} \Big[ \tr(X^x_u-X^{x'}_u)(b(X^x_u)-b(X^{x'}_u))du\\
\qquad \qquad \qquad + \tr(X^x_u-X^{x'}_u)\nabla\phi(X^x_u)dK_u^x-\tr(X^x_u-X^{x'}_u)\nabla\phi(X^{x'}_u)dK_u^{x'} \Big]\\
\qquad \qquad -\alpha e^{-\theta u}e^{-\alpha(\phi(X^x_u)+\phi(X^{x'}_u))} |X^x_u-X^{x'}_u|^2 \Big[ dK_u^x+dK_u^{x'}\\
\qquad \qquad \qquad + \tr(\nabla\phi(X^x_u)+\nabla\phi(X^{x'}_u))\sigma(d\tilde{W}_u+\beta_u du)\\
\qquad \qquad \qquad +\frac{1}{2}\textrm{Tr}(\nabla^2\phi(X^x_u)\sigma\tr\sigma+\nabla^2\phi(X^{x'}_u)\sigma\tr\sigma)du\\
\qquad \qquad  \qquad + \big(\tr\nabla\phi(X^x_u)b(X^x_u)+\tr\nabla\phi(X^{x'}_u)b(X^{x'}_u)\big)du \Big]\\
\qquad \qquad + \alpha^2 e^{-\theta u}e^{-\alpha(\phi(X^x_u)+\phi(X^{x'}_u))} |X^x_u-X^{x'}_u|^2 \Big[\\
\qquad \qquad \qquad \tr(\nabla\phi(X^x_u)+\nabla\phi(X^{x'}_u))\sigma\tr\sigma(\nabla\phi(X^x_u)+\nabla\phi(X^{x'}_u))\Big] ds.\\
\end{array}
\end{displaymath}
By Lemma~(\ref{lemme non convexe}) we have
$$ \Big(2\tr(X^x_u-X^{x'}_u)\nabla\phi(X^x_u) -\alpha |X^x_u-X^{x'}_u|^2 \Big)dK_u^x\leqslant 0,$$
and
$$ \Big(2\tr(X^{x'}_u-X^{x}_u)\nabla\phi(X^{x'}_u) -\alpha |X^x_u-X^{x'}_u|^2 \Big)dK_u^{x'}\leqslant 0.$$
Applying the definitions of $\beta$ and $\theta$, we obtain
\begin{displaymath}
\begin{array}{l}
d\Big(e^{-\theta u}e^{-\alpha(\phi(X^x_u)+\phi(X^{x'}_u))} |X^x_u-X^{x'}_u|^2\Big) \leqslant\\
-\alpha e^{-\alpha(\phi(X^x_u)+\phi(X^{x'}_u))} |X^x_u-X^{x'}_u|^2 \tr(\nabla\phi(X^x_u)+\nabla\phi(X^{x'}_u))\sigma d\tilde{W}_u.\\
\end{array}
\end{displaymath}
Thereby, for all $0 \leqslant t\leqslant s\leqslant n$
$$\mathbb{E}^{\mathbb{Q}_n}\Big[ e^{-\theta(s-t)-\alpha(\phi(X^x_s)+\phi(X^{x'}_s))} |X^x_s-X^{x^{\prime}}_s|\Big|\mathcal{F}_t\Big] \leqslant |X^x_t-X^{x^{\prime}}_t|.$$
The claim follows by setting $M=e^{2\alpha \sup_{x\in \bar{G}}\phi(x)}$. \cqfd
Of course we introduce a new hypothesis:
\paragraph{(H3'). } $\theta<0$.
\paragraph{}
\begin{thm}
Assume that $\sigma$ is a constant function. Theorems~\ref{existence convexe}, \ref{existence solution gdiff0}, \ref{existence solution gdiff0 mufixe} and \ref{existence solution gdiff0 mufixe 2} stay true if we substitute hypothesis (G2) and (H3) by (G2') and (H3').
\end{thm}

As in section~\ref{interpretation probabiliste de la solution d'une EDP}, it is possible to jointly modify $b$ and $\psi$ without modify the PDE \ref{EDP} if $\sigma$ is non-singular on $\overline{G}$. We set $\tilde{b}(x)=b(x)-\xi x$ and $\tilde{\psi}(x,z)=\psi(x,z)+\xi z \sigma^{-1} x$ for $\xi \in \mathbb{R}^+$. Then we are able to find a new hypothesis substituting (H3'). We note $\tilde{\theta}(\xi)$ the scalar $\theta$ corresponding to $\tilde{b}$ and $\tilde{\psi}$. Let $d$ the diameter of $\bar{G}$:
$$d:=\sup_{x,y \in \bar{G}} |x-y|.$$
\begin{prop}
$\tilde{\theta}(\xi) \leqslant \theta-(2-\frac{1}{2}d^2\alpha^2)\xi$. Particularly, if $\alpha d<2$ then there exists $\xi\geqslant0$ such that $\tilde{\theta}(\xi)<0$.
\end{prop}
\paragraph{Proof. }
Let $\tilde{\beta}$ the function $\beta$ linked with $\tilde{\psi}$. We have 
$$(Z^x_s-Z^{x'}_s)\tilde{\beta}_s=(Z^x_s-Z^{x'}_s)\beta_s+\frac{\xi}{2} (Z^x_s-Z^{x'}_s)\sigma^{-1}(X^{x'}_s+X^x_s)$$
So we can take $\tilde{\beta}_s=\beta_s+\frac{\xi}{2} \sigma^{-1}(X^{x'}_s+X^x_s)$. Thus $\tilde{\theta}(\xi) \leqslant \theta +C\xi$ with
\begin{eqnarray*}
C &=& -2+\sup_{x,y \in \bar{G}, x \neq y} \Big\{  -\dfrac{\alpha}{2}\tr(\nabla\phi(x)+\nabla\phi(y))(x+y)+\alpha(\tr\nabla\phi(x)x+\tr\nabla\phi(y)y)\Big\}\\
 &=& -2+\frac{\alpha}{2}\sup_{x,y \in \bar{G}} \big\{ \tr(\nabla\phi(x)-\nabla\phi(y))(x-y) \Big\}.
\end{eqnarray*}
On the other hand, we have
$$\sup_{x,y \in \bar{G}} \big\{ \tr(\nabla\phi(x)-\nabla\phi(y))(x-y) \big\}\leqslant d^2\alpha.$$
Indeed, according to the Taylor Lagrange theorem there exist $t,t' \in ]0,1[$ such that
$$\phi(x)=\phi(y)+\tr(x-y)\nabla\phi(y)+\frac{1}{2}\tr(x-y)\nabla^2\phi(ty+(1-t)(x-y))(x-y),$$
$$\phi(y)=\phi(x)+\tr(y-x)\nabla\phi(x)+\frac{1}{2}\tr(y-x)\nabla^2\phi(t'x+(1-t')(y-x))(y-x).$$
Finally $C \leqslant -2+\frac{d^2\alpha^2}{2}$ and the proof is therefore completed. \cqfd

\appendix
\section{Appendix}
\subsection{Proof of Proposition~\ref{convergence en loi de nu_n}}
We will prove that for all Lipschitz functions $f$, $\mathbb{E}^{\nu_n}[f] \stackrel{n \rightarrow +\infty}{\longrightarrow} \mathbb{E}^{\nu}[f]$. We set $X_0 \sim \nu$ and $X_0^n \sim \nu_n$, independent of $(W_t)_{t\geqslant 0}$. We have, for all $t\geqslant 0$,
\begin{equation*}
 \abs{\mathbb{E}^{\nu_n}[f]-\mathbb{E}^{\nu}[f]} \leqslant \underbrace{\abs{\mathbb{E}[f(X_t^{n,X_0^n})-f(X_t^{n,X_0})]}}_{A_{n,t}} + \underbrace{\abs{\mathbb{E}[f(X_t^{n,X_0})-f(X_t^{X_0})]}}_{B_{n,t}}.
\end{equation*}
Firstly,
$$A_{n,t} \leqslant K_f \mathbb{E} \abs{X_t^{n,X_0^n}-X_t^{n,X_0}}.$$
$\nabla^2 U_n \geqslant c I$, so $\nabla U_n$ is dissipative : we can prove that (see e.g. Proposition~3.3 in \cite{Fuhrman-Hu-Tessitore-07})
$$\mathbb{E} \abs{X_t^{n,X_0^n}-X_t^{n,X_0}} \leqslant e^{-ct} \mathbb{E} \abs{X_0^n-X_0}.$$
Then, by simple computations
\begin{eqnarray*}
 \mathbb{E} \abs{X_0^n-X_0} &\leqslant &\frac{1}{N}\int_{\mathbb{R}^d} \abs{x} e^{-U(x)} dx+ \mathbb{E} \abs{X_0} < + \infty.
\end{eqnarray*}
So, $A_{n,t} \leqslant Ce^{-ct} \stackrel{t \rightarrow +\infty}{\longrightarrow} 0$, and the limit is uniform in $n$. Moreover,
\begin{eqnarray*}
 B_{n,t} &\leqslant & K_f \mathbb{E} \abs{X_t^{n,X_0}-X_t^{X_0}} \leqslant K_f \int_{\overline{G}} \mathbb{E} [\sup_{0 \leqslant s \leqslant t} \abs{X_s^{n,x}-X_s^{x}}] \nu(dx).
\end{eqnarray*}
So, by Theorem~1 in \cite{Menaldi-Robin-85}, $B_{n,t} \stackrel{n \rightarrow +\infty}{\longrightarrow} 0$ when $t$ is fixed.  In conclusion, for all $t>0$,
$$\limsup_{n \rightarrow +\infty} \abs{\mathbb{E}^{\nu_n}[f]-\mathbb{E}^{\nu}[f]} \leqslant Ce^{-ct}.$$
So we can conclude the proof by letting $T \rightarrow +\infty$. \cqfd

\subsection{Proof of Proposition~\ref{generalisation ineg grandes dev}.}
We know that $\nabla^2 U_n \geqslant c I$. So, according to the Bakry-Emery criterion (see \cite{Bakry-Emery-85}), we have the Poincar\'e inequality
\begin{equation*}
\textrm{Var}_{\nu_n}(f)\leqslant -c^{-1}\langle \mathcal{L}_n f,f\rangle, \quad \forall f \in \mathbb{D}_2(\mathcal{L}_n).
\end{equation*}
Now, we are allowed to use Theorem~3.1 in \cite{Guillin-Leonard-Wu-Yao-07}:
\begin{equation*}
 \mathbb{P}\left(-\frac{1}{T}\int_0^T \mathcal{L}\phi(X_s^{n,X_0})ds \leqslant -\mathbb{E}^{\nu_n}[\mathcal{L}\phi]-\eps\right) \leqslant \mathbb{E}^{\nu}\left[\left(\frac{d\nu}{d\nu_n}\right)^2\right]^{1/2} \exp \left( -\frac{c\varepsilon^2 T}{\delta^2}  \right).
\end{equation*}
Firstly, by dominated convergence theorem
\begin{eqnarray*}
 \mathbb{E}^{\nu}\left[\left(\frac{d\nu}{d\nu_n}\right)^2\right]^{1/2} & = & \frac{N_n}{N} \stackrel{n \rightarrow +\infty}{\longrightarrow} 1.
\end{eqnarray*}
Moreover, applying Proposition~\ref{convergence en loi de nu_n},
$$ \mathbb{E}^{\nu_n}[\mathcal{L}\phi] \stackrel{n \rightarrow +\infty}{\longrightarrow} \mathbb{E}[\mathcal{L}\phi(X_0)].$$
Finally,
\begin{eqnarray*}
 \mathbb{E}\abs{\frac{1}{T}\int_0^T \mathcal{L}\phi(X_s^{n,X_0})ds-\frac{1}{T}\int_0^T \mathcal{L}\phi(X_s^{X_0})ds} & \leqslant & K_{\mathcal{L}\phi} \int_{\overline{G}}\mathbb{E}\left[\sup_{s \in [0,T]} \abs{X_s^{n,x}-X_s^x} \right] \nu(dx).
\end{eqnarray*}
But, according to \cite{Menaldi-Robin-85}, 
$$\mathbb{E}\left[\sup_{s \in [0,T]} \abs{X_s^{n,x}-X_s^x} \right] \stackrel{n \rightarrow +\infty}{\longrightarrow} 0$$
and the limit is uniform in $x$ belonging to $\overline{G}$. So
$$\mathbb{E}\abs{\frac{1}{T}\int_0^T \mathcal{L}\phi(X_s^{n,X_0})ds-\frac{1}{T}\int_0^T \mathcal{L}\phi(X_s^{X_0})ds} 
\stackrel{n \rightarrow +\infty}{\longrightarrow} 0,$$
and, as convergence in $L^1$ implies convergence in law, the claim follows. \cqfd


\end{document}